\documentclass[12 pt]{amsart}

\usepackage{amsmath, amssymb, amsthm}
\usepackage{tikz-cd}
\usepackage[shortlabels]{enumitem}
\usepackage[hidelinks]{hyperref}
\usepackage[pagewise]{lineno}

\newtheorem{theorem}{Theorem}[section]
\newtheorem{lemma}[theorem]{Lemma}
\newtheorem{corollary}[theorem]{Corollary}
\newtheorem{proposition}[theorem]{Proposition}
\newtheorem{remark}[theorem]{Remark}

\newtheorem{observation}[theorem]{Observation}

\newtheorem*{claim}{Claim}
\newtheorem{fact}[theorem]{Fact}

\newtheorem*{theorem*}{Theorem}
\newtheorem*{corollary*}{Corollary}

\theoremstyle{definition}
\newtheorem{example}[theorem]{Example}
\newtheorem{definition}[theorem]{Definition}

\newcommand{\abar}{\overline{a}}
\newcommand{\bbar}{\overline{b}}
\newcommand{\cbar}{\overline{c}}
\newcommand{\dbar}{\overline{d}}

\newcommand{\calL}{\mathcal{L}}

\newcommand{\calC}{\mathcal{C}}

\newcommand{\calU}{\mathcal{U}}

\newcommand{\wt}[1]{\widetilde{#1}}

\newcommand{\id}{\operatorname{id}}

\newcommand{\dl}{\mathrm{dlog}}
\newcommand{\logd}{\log_{\delta}}
\newcommand{\Aut}{\mathrm{Aut}}
\def\acl{\operatorname{acl}}
\def\dcl{\operatorname{dcl}}
\def\tp{\operatorname{tp}}

\newcommand{\bg}[3][\calC]{\Aut_{#2}(#3/#1)}

\def\forkindep{\mathrel{\raise0.2ex\hbox{\ooalign{\hidewidth$\vert$\hidewidth\cr\raise-0.9ex\hbox{$\smile$}}}}}

\newcommand{\partials}[2]{\frac{\partial #1}{\partial #2}}

\title[Internality of autonomous differential equations]{Internality of autonomous algebraic differential equations}

\author{Christine Eagles}
\address{Christine Eagles\\
University of Waterloo\\
Department of Pure Mathematics\\
Mathematics \& Computer\\
Waterloo, ON N2L 3G1\\
Canada}
\email{ceagles@uwaterloo.ca}

\author{L\'eo Jimenez}
\address{L\'eo Jimenez\\
The Ohio State University\\
Department of Mathematics\\
Math Tower\\
Columbus, OH 43210-1174\\
United States}
\email{jimenez.301@osu.edu}


\thanks{The first author was supported by the Natural Sciences and Engineering Research Council of Canada (NSERC) [reference number CGS D - 588616 - 2024]}

\date{\today}

\begin{document}

\begin{abstract}
This article is interested in internality to the constants of systems of autonomous algebraic ordinary differential equations. Roughly, this means determining when can all solutions of such a system be written as a rational function of finitely many fixed solutions (and their derivatives) and finitely many constants. If the system is a single order one equation, the answer was given in an old article of Rosenlicht. In the present work, we completely answer this question for a large class of systems. As a corollary, we obtain a necessary condition for generic solutions to be classical in the sense of Umemura. We then apply these results to determine exactly when solutions to Poizat equations (a special case of Li\'{e}nard equations) are internal, answering a question of Freitag, Jaoui, Marker and Nagloo, and to the classic Lotka-Volterra system, showing that its generic solutions are almost never classical (and thus not Liouvillian either).
\end{abstract}

\maketitle

\tableofcontents

\section{Introduction}

This article examines systems of order one ordinary differential equations:
\begin{align*}
\begin{cases}
    y_1' = f_1(y_1, \cdots , y_k) \\
    \vdots \\
    y_k' = f_k(y_1, \cdots , y_k) 
\end{cases} \tag{$\dagger$} \label{sys of eq intro}\end{align*}
\noindent where $f_1, \cdots , f_k$ are rational functions with coefficients in some differential field $F$ (all fields in this article will be of characteristic zero). Note that this is not as restrictive as it looks: many systems of higher order algebraic ordinary differential equations can be rewritten in this way. For example, a single equation of the form $(E): y^{(n+1)} = f(y^{(n)}, \cdots , y', y)$, for some rational function $f$, can be rewritten as the system:
\[
\begin{cases}
    y_0' = y_1 \\
    y_1' = y_2 \\
    \vdots \\
    y_{n}' = f(y_n, \cdots , y_1,y_0)
\end{cases}
\]
\noindent where the variable $y_i$ plays the role of the derivative $y^{(i)}$. For a more thorough examination of what systems can be transformed into a system of the form (\ref{sys of eq intro}), see Observation \ref{obs: gen-form} in the preliminaries section.

When the field $F$ is a field of constants, i.e. the derivative is identically zero on $F$, this system is autonomous, and it defines a rational vector field. We will almost always work in that context. 

We are interested in whether there are a finite number of solutions $a_1, \cdots , a_n$ and a rational function $h$ such that any generic solution $a$ of the system (\ref{sys of eq intro}) can be expressed as $a = h(c_1, \cdots , c_m, a_1 , \cdots , a_n)$, where the $c_i$ are complex numbers. This is an important phenomenon, which is exemplified by, but not limited to, systems of linear differential equations. It has been central in various areas of mathematics. 

Historically, we can cite Painlev\'{e}'s study of functions with \emph{integrals depending algebraically on constants} \cite{painleve1897leccons}, a similar definition, but replacing $h$ by an algebraic function. In model theory, it corresponds exactly to the notion of \emph{internality to the constants}, which we will define in the preliminaries. From a differential algebra perspective, it is closely related to Kolchin's \emph{strongly normal extensions} \cite{kolchin1953galois}. The connection between these two areas has been thoroughly investigated, see for example \cite[Chapter 2, Section 9]{marker2017model}. Moreover, this approach to differential equations has been investigated in depth by the Japanese school of differential algebra, for example Umemura in his study of \emph{classical functions} \cite{umemura1988irreducibility} and Nishioka, who clarified the relationship between Painlev\'{e}'s work and the perspectives of differential algebra and model theory in \cite{nishioka1989differential}. Finally, let us mention that similar ideas have been of interest to more applied branches of mathematics through the use of \emph{superposition laws} and \emph{Lie systems}, for example in physics \cite{carinena2008integrability}. 

There is also a connection to \emph{first integrals} of the associated vector field. Making the assumption that $F$ is a field of constants, if one considers the field $K = F(a_1, \cdots a_n)$ generated by the solutions used to express any other generic solution, there are elements $g_1, \cdots, g_k$ of $K(x_1, \cdots , x_k)$, algebraically independent over $K$, and such that their Lie derivative with respect to the associated vector field $X$, defined as
\[\calL_X(g_j) = \sum\limits_{i=1}^{k} \partials{g_j}{x_i}f_i \text{ ,}\]
is zero. These rational functions are called first integrals of the vector field associated to (\ref{sys of eq intro}), and we say that it has $k$ algebraically independent rational first integrals, after base change to $K$. This corresponds to a slight weakening of internality, called almost internality. 

Our question could thus be rephrased as follows: can we characterize rational vector fields of the form (\ref{sys of eq intro}) having $k$ algebraically independent rational first integrals after base change?

For the model-theorist, it is natural to work in a fixed differentially closed field to answer such questions, and this is the point of view we will take here. Note that any differential field embeds into a differentially closed field, and this is therefore an appropriate setting to study the existence of certain first integrals after base change. Moreover, by Seidenberg's embedding theorem \cite{seidenberg1958abstract}, any countably generated differential field embeds into a field of meromorphic function on some complex domain, which gives our work relevance for the study of complex analytic functions.

We will denote this fixed differentially closed field $\calU$ and let $\calC = \{ x \in \calU : \delta(x) = 0 \}$ be the \emph{constants} of $\calU$. It is an algebraically closed field. In the literature, the question of the internality of (\ref{sys of eq intro}) is only fully answered for a single order one autonomous differential equation in an old article of Rosenlicht \cite{rosenlicht1974nonminimality}:

\begin{theorem*}[Rosenlicht]
    Work in a saturated model $\calU \models \mathrm{DCF}_0$. Let $F$ be an algebraically closed field of constants and $f \in F(x)$. Then the generic type of $y' = f(y)$ is almost internal to the constants if and only if either $f = 0$, or $\partials{g}{x} f = 1$, or $\partials{g}{x}f = cg$ for some $g \in F(x)$ and $c \in F$.
\end{theorem*}

More recently, attempts have been made to obtain criteria for internality to the constants of higher order equations obtained as pullbacks of lower order ones (see \cite{jin2020internality}, \cite{jaoui2023relative} and \cite{eagles2024splitting}). However, these results have limited scope because they are conditional on expressing a differential equation as a pullback, as well as already knowing that the lower order equation is internal.

A proof of Rosenlicht's theorem can be obtained by using the differential Galois group, or binding group, associated to the equation. A geometric argument shows that this binding group is linear, and by failure of the inverse Galois problem for autonomous differential equations, must be commutative, and therefore either the additive or the multiplicative group. From this one obtains a definable map to the generic type of either $z' = 1$ or $z' = cz$, and an easy computation gives Rosenlicht's theorem. This is essentially our proof method for the general case, and as will become clear later, it specializes to order one to yield Rosenlicht's theorem.

Indeed, all these tools are available for systems of autonomous equations. Before stating our main theorem, let us sketch the ideas of the proof in order to make the statement more intuitive. We can show that the generic type of (\ref{sys of eq intro}), if internal to the constants, must have binding group either $G_m^l$ or $G_m^{l-1} \times G_a$. From this we can obtain a definable map to some system:
\begin{align*}\begin{cases}
    z_1' = \lambda_1 z_1 \\
    \vdots \\
    z_{l-1}' = \lambda_{l-1} z_{l-1} \\
    z_l ' = \lambda_l z \text{ or } z_l' = 1
\end{cases} \tag{$\bigstar$} \label{eq: linear}\end{align*}
\noindent where the last equation depends on the binding group and $l$ is some positive integer less or equal to $k$, the number of variables of the system (\ref{sys of eq intro}). 

Quantifier elimination in $\mathrm{DCF}_0$ tells us that the map is given by rational functions, which will be some of the $g_j$ of the theorem. If $l = k$ we are done, but if $l < k$ the dimension of the solution set of (\ref{eq: linear}) is strictly smaller than the dimension of the solution set of (\ref{sys of eq intro}), and we must fill the gap with rational first integrals, which will be the rest of the $g_j$. 

As stated before, the system (\ref{sys of eq intro}) defines a rational vector field $X$. For any $g \in F(x_1, \cdots , x_n)$, we can compute the Lie derivative of $g$ with respect to that vector field, which we recall is defined as:
\[\calL_X(g) = \sum\limits_{i=1}^{k} \partials{g}{x_i}f_i\]
\noindent and it is easy to see that if $(a_1, \cdots , a_k)$ is a solution of (\ref{sys of eq intro}), then $g(a_1, \cdots , a_n)' = \calL_X(g)(a_1, \cdots , a_n)$. Equipped with this, we can now state our main theorem:

\begin{theorem*}[In $\mathrm{DCF}_0$]
    Let $F$ be an algebraically closed field of constants, some $f_1, \cdots , f_k \in F(x_1, \cdots , x_k)$ and $p$ the generic type of the system:    
    \[\begin{cases}
    y_1' = f_1(y_1, \cdots , y_k) \\
    \vdots \\
    y_k' = f_k(y_1, \cdots , y_k)
    \end{cases}\]
    \noindent as well as $X$ the associated vector field. Then $p$ is almost internal to the constants if and only if there are $g_1 , \cdots , ,g_k \in F(x_1, \cdots , x_k)$, algebraically independent over $F$, such that for all $1 \leq i \leq k$, either:
    \begin{itemize}
        \item $\calL_X(g_i) = \lambda_i g_i$ for some $\lambda_i \in F$,
        \item $\calL_X(g_i) = 1$.
    \end{itemize}
    Moreover, the type $p$ is weakly $\calC$-orthogonal if and only if the $\lambda_i$ are $\mathbb{Q}$-linearly independent and the second equality is true for at most one $i$.
    
\end{theorem*}

From an integrability perspective, what we ask is the existence of a maximal number of algebraically independent rational maps such that their Lie derivatives satisfy very simple equations. We also take the opportunity to point out the connection with \emph{Darboux polynomials}: some $P \in F[x_1, \cdots , x_k]$ is a Darboux polynomial if there is another polynomial $Q \in F[x_1 , \cdots , x_k]$, called its \emph{cofactor}, such that $\mathcal{L}_X(P) = Q P$. The rational functions $g_i$ with $\mathcal{L}_X(g_i) = \lambda_i g_i$ satisfy a very similar equation: they are rational instead of polynomial, but must have constant cofactors. See the discussion in subsection 4.2 of \cite{freitag2023equations} for more details on Darboux polynomials.

Consider again the vector field interpretation of almost internality given earlier, i.e. that the vector field associated to (\ref{sys of eq intro}) has the maximum number of algebraically independent rational first integrals after base change. We can formulate a natural opposite condition: the associated vector field does not have any non-constant rational first integral, even after base change. This corresponds to the model-theoretic notion of \emph{orthogonality to the constants}. From a model-theoretic perspective, rational first integrals of the system correspond exactly to non-constant definable maps from its generic type to the constants, and we adopt that point of view in the next two paragraphs.

Internality of the system (\ref{sys of eq intro}) give rise to a very rich algebraic structure on its solution set. For example, it is acted upon by an algebraic group, which acts as a Galois group, which we use in the proof of our main theorem. Moreover, after base change, the system has a definable map to the constants (in model-theoretic language, this is after taking a non-forking extension), or equivalently, a rational first integral. 

On the other hand, orthogonality to the constants can be thought of as implying that solutions have a rather rudimentary structure. Indeed, the generic type of the system has no definable map to the constants, even after taking a non-forking extension. As a consequence of Zilber's dichotomy (proved by Hrushovski and Sokolovi\'{c} in \cite{hrushovski1994minimal}, see \cite{pillay2003jet} for an alternative proof by Pillay and Ziegler), it has no definable map to \emph{any} infinite field. However, there may be a definable map from a forking extension to the constants (an example of this is given by \cite[Example 5.6]{jaoui2023relative}). From a differential algebra perspective, it is known (see \cite[Theorem 4.6]{freitag2023equations}) that orthogonality to the constants implies that the generic solution is not classical in the sense of Umemura \cite{umemura1988irreducibility}, and in particular not Liouvillian.  

A strictly stronger non-structure property has recently been studied in the literature: strong minimality and triviality of the solution set. In this case, there is no definable map to the field of constants, even after taking a forking extension (indeed, the forking extension are all algebraic). Even more, no infinite group can be defined. Devilbiss and Freitag \cite{devilbiss2023generic} showed that the type of a generic solution to $f(x) = 0$, where $f$ is a differential polynomial of degree $d$ and order $h > 1$ with differentially transcendental independent coefficients, is strongly minimal as long as $d \geq 2(h+2)$ (triviality is still open, this is their Conjecture 1.7). Using the optimal bound on the nonminimality degree obtained by Freitag, Jaoui and Moosa \cite{freitag2023degree}, a better bound of $d \geq 6$ can be obtained. Jaoui proves strong minimality and triviality in \cite{jaoui2023density}, if the $f_i$ are polynomials with $\mathbb{Q}$-algebraically independent constant coefficients. However, this second result does not apply to many concrete systems of equations, as for example all coefficients must be non-zero. 

Our main theorem provides, as an easy corollary, a characterization of orthogonality to the constants for any autonomous system of the form (\ref{sys of eq intro}):

\begin{corollary*}\label{cor: intro}
   Let $F$ be an algebraically closed field of constants, some $f_1, \cdots f_k \in F(x_1, \cdots , x_k)$ and $p$ the generic type of the system:
    \[\begin{cases}
    y_1' = f_1(y_1, \cdots , y_k) \\
    \vdots \\
    y_k' = f_k(y_1, \cdots , y_k)
    \end{cases}\]
    \noindent as well as $X$ the associated vector field. Then $p$ is not orthogonal to the constants if and only if there is $g \in F(x_1, \cdots , x_k) \setminus F$ and some $\gamma \in F$ such that either $\calL_X(g) = \gamma g$ or $\calL_X(g)  = \gamma$. 
\end{corollary*}

Although this is a weaker conclusion than strong minimality and triviality, it can be readily applied to concrete systems of autonomous differential equations, and the rest of this article is mainly devoted to this. 

Before moving on, let us spell out a connection, pointed out to us by an anonymous referee, with a theorem of Jouanolou on algebraic foliations \cite{jouanolou1977hypersurfaces} and its (unpublished) extension by Hrushovski \cite{hrushovski1995odes}. See \cite{freitag2017finiteness} for a published proof and an extension to the partial differential case, as well as \cite{moosa2019model} for a nice survey of this and related results.

Let $a$ be a solution of $z' = \gamma z$ (respectively $z' = \gamma$), and $g$ be such that $\mathcal{L}_X(g) = \gamma g$ (respectively $\mathcal{L}_X(g) = \gamma$). Also let $P,Q \in F[x_1, \cdots , x_k]$ be relatively prime with $g = \frac{P}{Q}$. We can extend the Lie derivative $\mathcal{L}_X$ to $F(a)(x_1, \cdots , x_n)$ by letting $\mathcal{L}_X(a) = \gamma a$ (respectively $\mathcal{L}_X(a) = \gamma$). In both cases, it is straightforward to compute that $P-aQ$ divides $\mathcal{L}_X(P-aQ)$. We say that the hypersurface given by $P-aQ = 0$ is \emph{invariant}. Our corollary therefore states that if there are no invariant hypersurfaces of this form, then the generic type of (\ref{sys of eq intro}) is orthogonal to the constants.

In general, if $F<K$ is a differential field, there is a unique extension of $\mathcal{L}_X$ to $K(x_1, \cdots , x_k)$ coinciding with the derivative of $\mathcal{U}$ on $K$. If a hypersurface is given by the equation $S(x_1, \cdots x_k) = 0$, for some $S \in K[x_1, \cdots , x_k]$, we say it is invariant if $S$ divides $\mathcal{L}_X(S)$.

The previously mentioned theorems of Jouanolou and Hrushovski state, in particular, that orthogonality of (\ref{sys of eq intro}) to the constants implies that for any $F<K$, there are only finitely many invariant hypersurfaces. Therefore, our result here implies the following: \emph{if the system \emph{(\ref{sys of eq intro})} has no invariant hypersurfaces of the form given by the corollary, then over any $K$ containing $F$, it only has finitely many invariant hypersurfaces.}

We now summarize the article. Section \ref{sec: preliminaries} contains model theoretic preliminaries as well as reminders on logarithmic differential equations. Section \ref{sec: criteria} contains the proof of the main theorem. It is divided into two parts: the first, in Subsection \ref{subsec: weak orth}, we assume that the generic solution of (\ref{sys of eq intro}) is weakly orthogonal to the constants, i.e. has no rational first integral. This eliminates a number of complications, which are dealt with in Subsection \ref{subsec: not weak orth}, where the main theorem is finally proved. 

Section \ref{sec: applications} is devoted to applying the results of Section 
\ref{sec: criteria} to differential equations from the literature:
\begin{enumerate}
    \item We characterize exactly when generic solutions of Poizat equations
    \[y'' = y'f(y)\]
   with $f$ a rational function with constant coefficients, are almost internal to the constants. This answers a question of Freitag, Jaoui, Marker and Nagloo in \cite{freitag2023equations}.
    \item We show the generic solution of a Lotka-Volterra system
    \[
    \begin{cases}
        x' = ax + b x y \\
        y' = cy + d x y
    \end{cases}
    \]
    (with $a,b,c,d$ non-zero constants) is orthogonal to the constants, and in particular not classical, unless $a = c$. If $a = c$, the generic solution is elementary by a result of Varma \cite{varma1977exact}, and in model theoretic terms, is 2-analysable in the constants but not almost internal to the constants.

    Since this article was first written, Duan and Nagloo \cite{duan2025algebraic} proved that after discarding the degenerate solutions $x = 0$ and $y=0$, the set defined by this system of equations is strongly minimal whenever $\frac{a}{c} \not\in \mathbb{Q}$, which implies orthogonality to the constants. In a subsequent paper \cite{duan2025algebraicm}, Duan and the authors showed that it is also strongly minimal in the $a \neq c$ case. However, note that the computations and results of Subsection \ref{subsec: LV} are still required to classify the invariant algebraic curves in the $a = c$ case, see \cite[Theorem 3.13]{duan2025algebraicm}.
    \item We connect our methods to the pullbacks based criteria for internality developed in \cite{jin2020internality} and \cite{eagles2024splitting}.
\end{enumerate}

Finally in Section \ref{sec: generalizations}, we discuss two potential generalizations:

\begin{enumerate}[(a)]
    \item differential equations on arbitrary algebraic varieties,
    \item non-autonomous differential equations.
\end{enumerate}

In both cases, the main obstruction for our methods is Galois theoretic, as we loose our tight control over the Galois group. However, it seems likely that (a) could be studied using the Chevalley decomposition for algebraic groups. 

As for (b), we point out that a recent result of Jaoui and Moosa \cite{jaoui2022abelian} easily yields a generalization of Rosenlicht's Theorem to a single non-autonomous order one equation. Given an algebraically closed differential field $F$, it quickly follows from their work that the generic type of an equation $y' = f(y)$, for some $f \in F(x)$, is almost $\calC$-internal if and only if there is a non-constant $g \in F(x)$ and some $a,b,c \in F$ such that:
    \[f \partials{g}{x} + g^{\delta}= a g^2 + bg + c\]
\noindent where $g^{\delta}$ is the result of applying the unique derivation of $F(x)$ extending that of $F$, and with $\delta(x) = 0$, to the rational function $g$. This is Corollary \ref{cor: Jaoui-Moosa-cor} below. 

\section{Preliminaries}\label{sec: preliminaries}

\subsection{Model theory}

Throughout this article, we will work in a saturated model $(\calU, \delta)$ of the theory $\mathrm{DCF}_0$ of differentially closed fields, and we denote by $\calC$ its field of constants. For this preliminary section, fix some differential subfield $F< \calU$ which we will use as our base parameters. For any differential field $K < \calU$ and tuple $a$, we will denote $K(a)$ the field generated by $K$ and $a$, and $K \langle a \rangle$ the differential field generated by $K$ and $a$.

We will have to deal with systems of differential equations in $\calU$ as well as rational functions and Laurent series fields. We will save the variables $y,z,u,v$ for differential equations in $\mathrm{DCF}_0$, and the variables $x, x_0,x_1,x_2 , \cdots $ will be used for rational functions and Laurent series.

We assume familiarity with geometric stability theory, for which a good reference is \cite{pillay1996geometric}. In particular, chapter 7 of that book contains most of what we need regarding the key notion of internality. We also assume familiarity with the interaction of model theory and differentially closed fields, for which a good reference is \cite[Chapter 2]{marker2017model} (see also the preliminaries of \cite{freitag2023equations} for a succinct introduction).

By definable, we always mean with parameters, and will always specify the parameters when relevant. When $p$ is a type, we will denote $p(\calU)$ its set of realizations in $\calU$. Often, we will identify a type and its set of realizations in $\calU$ and say, for example, that a group $G$ acts on a type $p$ to mean that $G$ acts on $p(\calU)$. We will also, given two types $p$ and $q$, talk about definable maps $f : p \rightarrow q$. By that we mean a partial definable function such that $\mathrm{dom}(f) \supset p(\calU)$ and $f(p(\calU)) \subset q(\calU)$ (if $f,p$ and $q$ are over the same parameters, this implies $f(p(\calU)) = q(\calU)$). Similarly, we may write about definable maps $f : p \rightarrow \calC$.

We are interested in obtaining necessary and sufficient conditions for certain types to be \emph{internal to the constants}:

\begin{definition}
    The type $p \in S(F)$ is said to be internal to the constants, or $\calC$-internal, if it is stationary and there is an extension $F<K$ such that for any $a \models p$, there is a tuple of constants $c_1, \cdots , c_n$ such that $a \in \dcl(c_1 \cdots c_n K)$. 

    If we replace $\dcl$ by $\acl$, the type $p$ is said to be \emph{almost} $\calC$-internal.
\end{definition}

If $p$ is $\mathcal{C}$-internal then the extension $K$ can be chosen to be $F\langle a_1,...,a_n\rangle$, where $a_1, \cdots , a_n$ is an $F$-independent tuple of realizations of $p$. We call such a tuple a fundamental system of solutions for $p$. We say that $p$ is fundamental if some (equivalently any) $a \models p$ is a fundamental system of solutions for $p$.

We will use the following well-known fact, which is not quite immediate from the definition:

\begin{lemma}\label{lem: int-bir-to-var}
    Let $p \in S(F)$ be some $\calC$-internal type. Then there is some $F<K$, some $l \in \mathbb{N}$, some $K$-type definable $X \subset \calC^l$, and a $K$-definable bijection $f: X \rightarrow p$.
\end{lemma}

\begin{proof}
    By \cite[Chapter 7, Lemma 4.2]{pillay1996geometric}, there are $a_1, \cdots , a_n \models p$ and an $a_1 \cdots a_n F$-definable partial function $f$ such that for any $a \models p$, there are $c_1, \cdots , c_m \in \calC$ such that $f(c_1, \cdots, c_m) = a$.
    
    Let $K = F\langle a_1, \cdots , a_n \rangle$ and $Y = \{ (c_1, \cdots , c_n) \in  \calC^n, f(c_1, \cdots , c_n) \models p \}$. We define an equivalence relation $E$ by $E(\cbar,\dbar)$ if $f(\cbar) = f(\dbar)$ (or $f$ is not defined at $\cbar$ and $\dbar$). Then $X = Y/E$ is a $K$-type definable set by elimination of imaginaries, and by construction, there is a $K$-definable bijection between $p$ and $X$.

\end{proof}

Essential to us will be the classical fact that internal types give rise to definable groups of automorphisms. More precisely, let $\bg{F}{\calU}$ be the group of autmorphisms of $\calU$ fixing $F \cup \calC$ pointwise. Given any type $p \in S(F)$, we can consider the following group of automorphisms:
\[\bg{F}{p} := \left\{ \sigma \vert_{p(\calU)}: \sigma \in \bg{F}{\calU} \right\}\]
\noindent i.e. restrictions to $p(\calU)$ of automorphisms fixing $F \cup \calC$ pointwise.

It is called the \emph{binding group of $p$ over $\calC$}, and we have:

\begin{fact}\label{fact: def-binding}
    Let $p \in S(F)$ be a $\calC$-internal type. Then the action of $\bg{F}{p}$  on $p$ is isomorphic to a canonical $F$-definable group action.
\end{fact}

See \cite[Chapter 7, Theorem 4.8]{pillay1996geometric} for a proof of this. In the rest of this article, when we write $\bg{F}{p}$, we will mean the definable group of Fact \ref{fact: def-binding}. As a consequence of \cite[Chapter 7, Remark 4.9]{pillay1996geometric} applied to this context, it is definably isomorphic to the $\mathcal{C}$-points of an algebraic group:

\begin{fact}\label{fact: int-equiv-map}
    Let $p \in S(F)$ be a $\calC$-internal type. Then its binding group $\bg{F}{p}$ is definably (maybe using some extra parameters) isomorphic to $G(\calC)$, the $\calC$-points of some algebraic group $G$. In particular $G(\calC)$ acts definably on $p$.
\end{fact}

We will sometimes make a slight abuse of terminology and talk about properties of the algebraic group $G$ as if they were properties of $\bg{F}{p}$. For example, we may say that $\bg{F}{p}$ is linear to mean that is it definably isomorphic to the $\calC$-points of a linear algebraic group. This is harmless because as a consequence of stable embeddedness of $\calC$, if $\bg{F}{p}$ is isomorphic to $G(\calC)$ and $H(\calC)$ for algebraic groups $G$ and $H$, then $G$ and $H$ are isomorphic \emph{as algebraic groups}.

A type is internal to the constants when it is, after extending the parameters, in definable bijection with the type of a tuple of constants. We could ask about opposite properties: the type $p$ does not have any definable relationship with constants (with or without allowing extra parameters). This intuition is made precise by the following:

\begin{definition}
    A type $p \in S(F)$ is:
    \begin{itemize}
        \item weakly $\calC$-orthogonal (or weakly orthogonal to $\calC$) if any realization $a$ of $p$ is independent, over $A$, from any tuple of constants,
        \item orthogonal to $\calC$ if is is stationary and for any $B \supset A$, any $a \models p\vert_B$ is independent, over $B$, from any tuple of constants. 
    \end{itemize}
\end{definition}

Even though our starting motivation was internality, our methods will also give criteria for weak $\calC$-orthogonality and $\calC$-orthogonality. A non-trivial fact, on which we will expand in Lemma \ref{lem: weak-ortho-cored}, is that $p \in S(F)$ is not weakly $\calC$-orthogonal if and only if there is an $F$-definable function $f : p \rightarrow \calC$ with $f(a) \not\in F^{\mathrm{alg}}$ for any $a \models p$ (see \cite[Lemma 2.1]{freitag2023bounding} for a proof). In other words, non weak $\calC$-orthogonality will correspond to the existence of a first integral given by a rational function of a solution and its derivatives. From this, it is straightforward to deduce a concrete criteria for weak orthogonality of certain types, and we do so in Proposition \ref{prop: weak-ortho-crit}.

No non-algebraic $\calC$-internal type can be $\calC$-orthogonal. However it could be weakly $\calC$-orthogonal, and this is characterized by:

\begin{fact}
    The $\calC$-internal type $p \in S(F)$ is weakly $\calC$-orthogonal if and only if $\bg{F}{p}$ acts transitively on $p$.
\end{fact}

See \cite[Lemma 2.15]{jin2019logarithmic} for a proof of this. Key to our methods will be a strong restriction on $G$ when the binding group is linear and $F$ is an algebraically closed field of constants (see the proof of \cite[Theorem 3.9]{freitag2022any}):

\begin{fact}\label{fact: if bind-grp lin then GmGa}
    Assume $F = F^{\mathrm{alg}} < \calC$. Let $p \in S(F)$ be a $\calC$-internal and weakly $\calC$-orthogonal type. If the binding group of $p$ is definably isomorphic to the $\calC$-points of a linear algebraic group, then it is definably isomorphic to the $\calC$-points of $(G_m)^l \times (G_a)^k$ for some $l \in \mathbb{N}$ and $k \in \{ 0 , 1 \}$.
\end{fact}

As mentioned in the introduction, we will be interested in systems of equations of the form:
\begin{align*}\begin{cases}
    y_1' = f_1(y_1, \cdots , y_k) \\
    \vdots \\
    y_k' = f_k(y_1, \cdots , y_k) 
\end{cases} \tag{$\dagger$} \label{sys of eq prel}\end{align*}
\noindent for some $f_1, \cdots , f_k \in F(x_1, \cdots, x_k)$. Many systems can be transformed into one of the form (\ref{sys of eq intro}):

\begin{observation}\label{obs: gen-form}
    Let $f_1, \cdots , f_k$ be rational functions over $F$ and $n_1, \cdots , n_k$ be positive integers. The solutions of
    \[
    \begin{cases}
        y_1^{(n_1+1)} = f_1(y_1^{(n_1)}, \cdots ,y_1, \cdots y_k^{(n_k)}, \cdots, y_k) \\
        \vdots \\
        y_k^{(n_k+1)} = f_k(y_1^{(n_1)}, \cdots , y_1, \cdots y_k^{(n_k)}, \cdots, y_k)
    \end{cases}
    \]
    are in bijection with solutions of the system of equations:
    \[
    \begin{cases}
        y_{i,j}' = y_{i,j+1}, 1 \leq i \leq k, 0 \leq j \leq n_i \\
        y_{i,n_i}' = f_i(y_{1,n_1}, \cdots y_{1,0}, \cdots y_{k,n_k}, \cdots, y_{k,0}), 1 \leq i \leq k
    \end{cases}
    \]
    through the function sending $y_i^{(j)}$ to $y_{i,j}$ for all $1 \leq i \leq k$ and all $0 \leq j \leq n_i$.
\end{observation}

As a simple example of such a transformation, see Subsection \ref{subsec: poizat} on Poizat equations. 

We would rather work with types than definable sets, and we will therefore always consider the \emph{generic type} of such a system: it is the type $p \in S(F)$ of a solution $a$ of (\ref{sys of eq prel}) that does not satisfy any polynomial equation over $F$. In other word, such that $F \langle a \rangle$ is isomorphic, over $F$, to the differential field $F(x_1, \cdots , x_k)$ with $\delta(x_i) = f_i(x_1, \cdots , x_k)$.

We will show that if it is $\calC$-internal, the generic type $p$ of such a system is \emph{interdefinable} with very specific types:

\begin{definition}
    Let $p,q \in S(F)$, we say that $p$ and $q$ are interdefinable if for any $a \models p$, there is $b \models q$ such that $\dcl(aF) = \dcl(bF)$. This implies that there is an $F$-definable bijection between $p(\calU)$ and $q(\calU)$.
\end{definition}

By quantifier elimination for $\mathrm{DCF}_0$, an $F$-definable function with domain the generic type $p$ of a system of the form (\ref{sys of eq prel}) is given by rational functions in the $y_i$. We will make frequent use of this.

Finally, the more subtle notion of \emph{analysability in the constants} will make an appearance:

\begin{definition}
    A type $p \in S(F)$ is analysable in the constants, or $\calC$-analysable, if there is $a_1 = a \models p, a_2, \cdots , a_n$ such that:
    \begin{itemize}
        \item for all $1 \leq i \leq n-1$, the type $\tp(a_i/a_{i+1}F)$ is $\calC$-internal,
        \item for all $1 < i \leq n$ we have $a_i \in \dcl(a_{i-1}F)$,
        \item the type $\tp(a_n/F)$ is $\calC$-internal.
    \end{itemize}

    This gives rise to a sequence of types and definable functions (all over $F$):
    \[p = p_1 \xrightarrow{f_1} p_2 \xrightarrow{f_2} \cdots p_{n-1} \xrightarrow{f_{n-1}} p_n\]
    \noindent where $p_i = \tp(a_i/F)$.
    We call $n$ the number of \emph{steps} of the analysis, and say that $p$ is $\calC$-analysable in $n$ steps, or $n$-analysable over the constants.
\end{definition}

In the above definition, we could also replace internality to the constants with internality to any minimal type. In this case, we obtain the definition of a \emph{semi-minimal analysis}. This will not play an important role in our work and we therefore do not say more about this.

Note that a type $p \in S(F)$ is $\calC$-internal if and only if $p$ is $\calC$-analysable in one step. But even for $n = 2$, there are many types that are $\calC$-analysable but not (almost) $\calC$-internal. Finding and characterizing such types has been an important theme in model theoretic differential algebra, see for example \cite{jin2018constructing}, \cite{jin2020internality} and \cite{eagles2024splitting}. Our methods allow us to give new examples of such types: a family of Weierstrass equations in Example \ref{ex: weierstrass} and some special case of Lotka-Volterra systems in Theorem \ref{theo: LV-is-ortho}.

\subsection{Logarithmic-differential equations}

We now recall some results on tangent bundles and derivations, as well as the definition of the logarithmic derivative.

Any affine algebraic variety $V \subset \calU^n$ has a tangent bundle $TV$, given by:
\[\Bigl\{ (a,u) \in \calU^n:  a \in V, \sum\limits_{i = 1}^n \frac{\partial P}{\partial x_i}(a)u_i = 0 \text{ for } P \in I(V)  \Bigr\} \text{ .}\]
We have a natural projection $\pi : TV \rightarrow V$, and for any $a \in V$, we let $T_a V = \pi^{-1}(\{ a \})$. If $V$ is over a field of constants, the derivation $\delta$ gives a section of the tangent bundle by:
\begin{align*}
    \nabla : V & \rightarrow T(V)\\
    a & \rightarrow (a, \delta(a))
\end{align*}
By patching, we obtain a tangent bundle and $\nabla$ map for any algebraic variety. 

Given a regular map $f : V \rightarrow W$ and $a \in V$, there is a linear map $d f_a : T_a V  \rightarrow T_{f(a)} V $ given by differentiating $f$. We thus obtain a map $df : TV \rightarrow TW$. This gives a functor on the category of algebraic varieties. 

In particular, suppose that $G$ is a connected algebraic group, we consider, for $g \in G$, the map $\lambda^g(x) = gx$. We have a map:
\begin{align*}
    \eta : TG & \rightarrow T_e G \\
    (g,u) & \rightarrow (d\lambda^{g^{-1}})_g (g,u)
\end{align*}
If $G$ is defined over a field of constants, the \emph{logarithmic derivative} on $G$ is given by $\eta \circ \nabla = \dl_G$. A \emph{logarithmic differential equation} on $G$ is an equation of the form $\dl_G(x) = (e,u)$ for some $(e,u) \in T_e G(\calU)$. The set of solution to any such equation is $\calC$-internal, and thus so is its generic type $p$. The logarithmic differential equation is said to be \emph{full} if in addition $p$ is weakly $\calC$-orthogonal. In that case, it is isolated by $\dl_G(x) = (e,u)$ (and $x \in G$) and its binding group is definably isomorphic to $G(\calC)$ (in general, it is always isomorphic to a subgroup of $G(\calC)$).

Logarithmic differential equations will be key to our methods because of the following fact, which is a model-theoretic translation, given by Jaoui and Moosa in \cite{jaoui2022abelian}, of a theorem of Kolchin:

\begin{fact}(Kolchin)\label{fact: kolchin log-diff}
    Let $F$ be an algebraically closed differential field, and $p \in S(F)$ a $\calC$-internal, fundamental  and weakly $\calC$-orthogonal type. Then there is a connected algebraic group $G$, defined over $\calC \cap F$, such that $p$ is interdefinable with the generic type of a full logarithmic-differential equation on $G$ over $F$.

\end{fact}

We will only be interested in logarithmic differential equations on very basic algebraic groups, namely groups of the form $G = G_a ^m \times G_m ^n$ for some $n,m \in \mathbb{Z}^+$. In that case, first note that $T( G_a ^m \times G_m ^n )= (T G_a)^m \times (T G_m)^n$, and if $(e,u) = (e, u_1, \cdots , u_m, v_1, \cdots , v_m) \in T_e G$, then the equation $\dl_G(x) = (e,u)$, with $x = (x_1, \cdots , x_m, y_1, \cdots  , y_n)$ becomes a system of equations (we will now ignore the $e$ coordinate):
\[\begin{cases}
    \dl_{G_a}(x_i) = u_i , 1 \leq i \leq m\\
    \dl_{G_m}(y_j) = v_j , 1 \leq j \leq n \text{ .}
\end{cases}\]
\noindent Moreover it is easy to compute that $\dl_{G_a}(x) = \delta(x)$ for any $x \in G_a$, and $\dl_{G_m}(x) = \frac{\delta(x)}{x} : = \logd(x)$ for any $x \in G_m$ ($\logd$ is often called the logarithmic derivative in the literature).

The logarithmic differential equation $\dl_G (x) = u$ is full exactly when this system of equations has binding group $G_a^m \times G_m ^n$, and we need to characterize when this is the case. We first deal with the case of a product of $G_m$. 

\begin{lemma}\label{lem: fullness for Gm}
    Let $F < \calC$ be an algebraically closed field and $v_1, \cdots , v_n \in F$. The binding group of the generic type $q \in S(F)$ of the system of equations:
    \[\Bigl\{ \dl_{G_m}(y_j) = v_j , 1 \leq j \leq n \Bigr\}\]
    is isomorphic to $G_m(\calC)^n$ if and only if the $v_j$ are $\mathbb{Q}$-linearly independent.
\end{lemma}

\begin{proof}
    Suppose first that the $v_j$ are $\mathbb{Q}$-linearly independent. Let $\abar \models q$, because $q$ is generic, we have that $F\langle \abar \rangle = F(\abar)$ is isomorphic, over $F$, to the differential field $F(x_1, \cdots , x_n)$ with $\delta(x_j) = v_j x_j$ for all $j$. By \cite[Example 3.23]{magid1994lectures}, the field $F(x_1, \cdots , x_n)$ is a Picard-Vessiot extension of $F$. This implies, in particular, that it has no new constants and thus that $q$ is weakly $\calC$-orthogonal. By the previous discussion, its binding group must be definably isomorphic to $G_m(\calC)^n$.

    Conversely, suppose that the $v_j$ are not $\mathbb{Q}$-linearly independent. Let $q_j$ be the generic type of $\logd (y_j) = v_j$. There are $\lambda_j \in \mathbb{Z}$ such that $\sum\limits_{j= 1}^n \lambda_j v_j = 0$. Then one can compute that, for any $a_1 , \cdots , a_n$ with $\logd(a_j) = v_j$, we have that $\prod\limits_{j = 1}^n a_j^{\lambda_j} \in \calC$, which implies that $a_n \in \acl(a_1, \cdots , a_{n-1}, \calC)$, and thus $\bg{F}{q}$ is a finite extension of $ \bg{F}{q_1 \otimes \cdots \otimes q_{n-1}}$, which cannot be definably isomorphic to $G_m(\calC)^n$ as $ \bg{F}{q_1 \otimes \cdots \otimes q_{n-1}}$ is definably isomorphic to a definable subgroup of $G_m(\calC)^{n-1}$. 
\end{proof}

We now prove a version for $G_a$:

\begin{lemma}\label{lem: fullness for Ga}
    Let $F < \calC$ be an algebraically closed field and $u_1, \cdots , u_n \in F$. The binding group of the generic type $q \in S(F)$ of the system of equations:
    \[\Bigl\{ \dl_{G_a}(y_j) = u_j , 1 \leq j \leq m \Bigr\}\]
    is isomorphic to $G_a(\calC)^m$ if and only if $m = 1$ and $u_1 \neq 0$. 

    Moreover, in that case, if $F$ is a field of constant, then $q$ is interdefinable with the generic type of $\delta(y) = 1$.
\end{lemma} 

\begin{proof}
    Since $\bg{F}{q}$ must be an isomorphic to a definable subgroup of $G_a(\calC)^m$, it is well-known, for example by the proof of \cite[Theorem 3.9]{freitag2022any}, that it must in fact be isomorphic to a definable subgroup of $G_a(\calC)$. Thus it is either isomorphic to $G_a(\calC)$ or trivial. This already forces $m =1$. 

    So all that is left to show is that the generic type $q$ of $\delta(y) = u$, for $u \in \calC$, has binding group $G_a$ is and only if $u \neq 0$. If $u = 0$, then $q$ has only constant realizations, so the binding group is trivial. If the binding group is not $G_a$, it must be trivial, and thus for any $a \models q$, we have $a \in \dcl(\calC) = \calC$. Thus $u = \delta(a) = 0$. 

    For the moreover part, it is immediate that $y \rightarrow \frac{y}{u}$ is an $F$-definable bijection between the generic type of $\delta(y) = u $ and $\delta(y) = 1$. 
\end{proof}

Finally, we can combine these two lemmas into:

\begin{lemma}\label{lem: fullness}
    Let $G = G_a ^m \times G_m ^ n$, and consider a logarithmic differential equation on $G$ over an algebraically closed $F < \calC$, which we write as a system:
    \[\begin{cases}
    \dl_{G_a}(x_i) = u_i , 1 \leq i \leq m\\
    \dl_{G_m}(y_j) = v_j , 1 \leq j \leq n
    \end{cases}\]
    \noindent Then this equation is full if and only if both:
    \begin{itemize}
        \item $m \leq 1$, and if $m=1$ then $u_1 \neq 0$,
        \item $v_1, \cdots v_n$ are $\mathbb{Q}$-linearly independent.
    \end{itemize}
\end{lemma}

\begin{proof}
    First suppose that the equation is full. Then both the $\dl_{G_a}$ and $\dl_{G_m}$ systems must be full, and we conclude using Lemmas \ref{lem: fullness for Gm} and \ref{lem: fullness for Ga}.

    Conversely, suppose we have the two listed conditions. Then by Lemmas \ref{lem: fullness for Gm} and \ref{lem: fullness for Ga}, the $\dl_{G_a}$ and $\dl_{G_m}$ subsystems must be full, i.e. their generic types have binding groups isomorphic to $G_a(\calC)$ and $G_m(\calC)^n$, respectively. The generic type of the whole system has binding group definably isomorphic to a definable subgroup $G < G_a(\calC) \times G_m(\calC)^n$, and by fullness of the two subsystems, the group $G$ projects surjectively onto both $G_a(\calC)$ and $G_m(\calC)^n$. This implies that $G = G_a(\calC) \times G_m(\calC)^n$.
\end{proof}

\section{Criteria for internality}\label{sec: criteria}

\subsection{Weakly orthogonal case}\label{subsec: weak orth}

For the rest of this section, we will fix some algebraically closed field of constants $F$ and consider systems of differential equations of the form:
\begin{align*}\begin{cases}
    y_1' = f_1(y_1, \cdots , y_k) \\
    \vdots \\
    y_k' = f_k(y_1, \cdots , y_k)
\end{cases} \tag{$\dagger$} \label{sys of eq}\end{align*}
\noindent where $f_i \in F(x_1, \cdots , x_k)$ for all $i$.

We first prove a lemma to show the binding groups we care about are linear:

\begin{lemma}\label{lem: lin-bind-group-im}
    Let $p \in S(F)$ be the generic type of a system of the form \emph{(\ref{sys of eq})} and $\pi : p \rightarrow r$ be an $F$-definable map. If $r$ is $\calC$-internal, then its binding group is linear. 
\end{lemma}

\begin{proof}

    By Fact \ref{fact: int-equiv-map}, there is a definable bijection (maybe over extra parameters) between $r$ and some type-definable set of constants $X$. Moreover, the binding group $\bg{F}{r}$ is definably isomorphic, over extra parameters, to $G(\calC)$, where $G$ is an algebraic group. Let $F< K$ over which both of these maps are defined. 

    Because $X$ is in definable bijection with a type, it must be given by a Zariski dense subset of a variety $V$ defined over $K$. We can transfer the action of $G(\calC)$ on $p$ to an action of $G(\calC)$ on $X$.
    In this situation, Weil's regularization theorem, also known as Weil's group chunk theorem, states that there is an algebraic variety $W$ defined over $K$ on which $G$ acts regularly, and a $G$-equivariant birational map between $V$ and $W$. See Theorem 4.11 of \cite{bouscaren2009model} for more details. Because the action of $G$ on $X$ is faithful, so is the action of $G$ on $W$.

    Consider the \emph{Albanese morphism} $W \rightarrow A(W)$, which is a morphism to the abelian variety $A(W)$, characterized by the universal property that any morphism from $W$ to an abelian variety factors through it. By a result of Matsumura \cite{matsumura1963algebraic} (see also \cite{brion2010some}), the action of $G$ on $W$ induces a morphism $G \rightarrow A(W)$, and its kernel is linear. 

    If we show that $W$ is unirational, i.e. that its function field is a subfield of a purely transcendental expansion of $K$, then we could conclude, by \cite[Theorem 4]{lang2019abelian} and the corollary following it, that the map $W \rightarrow A(W)$ is constant, i.e. $A(W)$ is trivial. This implies in particular that the kernel of $G \rightarrow A(W)$ is $G$ itself, and thus $G$ is linear. We now prove that $W$ is unirational.

    Fix some $b \models r\vert_K$. There are $\abar \models p\vert_K$ and $g_1, \cdots , g_n \in F(x_1, \cdots, x_k)$ such that $b = \pi(\abar) = (g_1(\abar), \cdots , g_n(\abar))$. We know that $U(\abar/F)$ is finite, therefore $U(b/F)$ is finite as well, and in particular there is $l \in \mathbb{N}$ such that $F\langle b \rangle < F(b , b', \cdots , b^{(l-1)})^{\mathrm{alg}}$. 

    Let $P$ be the minimal polynomial of $b^{(l)}$ over $F(b, b' , \cdots , b^{(l-1)})$. By differentiating the equation $P(b^{(l)}) = 0$, we obtain that there is a rational function $h \in F(x_0, \cdots , x_{l})$ such that $b^{(l+1)} = h(b, b' , \cdots , b^{(l)})$. In particular:
    \begin{align*}
        K\langle b \rangle & = K ( b, b ' , \cdots , b^{(l+1)} ) \\
        & = K(\pi(\abar), \pi(\abar)', \cdots , \pi(\abar)^{(l+1)}) \\
        & < K(\abar) \text{ .}
    \end{align*}

    The field $K(\abar)$, since $a \models p\vert_K$, is a purely transcendental extension of $K$. But the field $K\langle b \rangle$ is isomorphic to the function field of $W$ over $K$, and therefore the latter is a subfield of a purely transcendental extension of $K$.

\end{proof}

We note the following consequence:

\begin{corollary}\label{cor: non-linear implies non-int}
    Let $p$ be the generic type of a system of the form \emph{(\ref{sys of eq})}. If there is a semi-minimal analysis $p = p_1 \xrightarrow{f_1} p_2 \xrightarrow{f_2}  \cdots \xrightarrow{f_n} p_{n+1} $ of $p$ such that for some $i \in \{ 1, \cdots , n+1 \}$ and $a \models p_i$ the binding group $\bg{f_i(a) F}{\tp(a/f_i(a)F)}$ is isomorphic to a non-linear algebraic group, then $p$ is not $\calC$-internal.
\end{corollary}

\begin{proof}
    
    Suppose, for a contradiction, that $p$ is $\calC$-internal. Then by Lemma \ref{lem: lin-bind-group-im}, its binding group must be linear, let $G$ be the linear algebraic group such that $\bg{F}{p}$ is definably isomorphic to $G(\calC)$. 

    Let $\pi $ be defined as the composition $f_{i-1} \circ \cdots \circ f_1$, it is an $F$-definable map from $p$ to $p_i = \pi(p)$. It induces a definable group morphism $\wt{\pi}$ fitting in a short exact sequence:
    \[1 \rightarrow \ker(\wt{\pi}) \rightarrow \bg{F}{p} \xrightarrow{\wt{\pi}} \bg{F}{\pi(p)} \rightarrow 1 \text{ .}\]
    As Lemma \ref{lem: lin-bind-group-im} tells us that the binding group of $\pi(p)$ is still linear, we replace $p$ by $\pi(p)$, i.e. we assume that $p$ is $\calC$-internal with an $F$-definable fibration $f: p \rightarrow f(p)$ such that $\bg{f(a)F}{\tp(a/f(a)F)}$ has binding group definably isomomorphic to the $\calC$-point of a non-linear algebraic group.

    We again get an induced $F$-definable morphism $\wt{f} : \bg{F}{p} \rightarrow \bg{F}{f(p)}$, and by \cite[Lemma 2.7]{jaoui2022abelian} and our assumption on the group $\bg{f(a)F}{\tp(a/f(a)F)}$, its kernel is $F$-definably isomorphic to the $\calC$-points of a non-linear algebraic group. This is a contradiction, as it must be definably isomorphic to a definable subgroup of $G(\calC)$ and $G$ is linear.

\end{proof}

This can be used to give examples of types that are 2-analysable in the constants but not internal to the constants:

\begin{example}\label{ex: weierstrass}
    Let $p$ be the generic type of $y'' = 6y^2 - \frac{1}{2}$, which is interdefinable with the generic type of the system:
    \[\begin{cases}
    y_0' = y_1\\
    y_1' = 6y_0^2 - \frac{1}{2}
\end{cases}\]

\noindent We define the function $f(y) = (y')^2-4y^3+y$, an easy computation shows that for any $a \models p$, we have that $f(a) \in \calC$. The fiber above $f(a)$ is the generic type of the Weiestrassian equation $(y')^2 = 4y^3 - y + f(a)$, well-known to be $\calC$-internal with binding group an elliptic curve (see \cite[Chapter 2, Section 9]{marker2017model}, for example). The map $f$ gives rise to an analysis of $p$ over $\calC$ in two steps. By Corollary \ref{cor: non-linear implies non-int}, the type $p$ cannot be $\calC$-internal.
\end{example}

Of crucial importance to us is the following immediate corollary of Facts \ref{fact: if bind-grp lin then GmGa} and \ref{fact: kolchin log-diff}:

\begin{corollary}\label{cor: interdef-full-log}
    Let $F$ be an algebraically closed field of constants. If $p$ is $\calC$-internal, weakly $\calC$-orthogonal and $\bg{F}{p}$ is definably isomorphic to the $\calC$-points of a linear algebraic group, then $p$ is fundamental and interdefinable with the generic type of a full logarithmic-differential equation over $F$ on either $G_a \times (G_m)^{k-1}$ or $(G_m)^k$, where $k = U(p)$.
\end{corollary}

\begin{proof}
    By Fact \ref{fact: if bind-grp lin then GmGa}, the binding group of $p$ is abelian, and either $G_m^k$ or $G_m^{k-1} \times G_a$ for some $k$. Note that the binding group action is always faithful. But a faithful transitive action of an abelian group is always regular, and thus $p$ is fundamental.

    By Fact \ref{fact: kolchin log-diff}, there is connected algebraic group $H$, defined over $\calC \cap F$, such that $p$ is interdefinable with the generic type of a full logarithmic-differential equation on $H$ over $F$. This forces $\bg{F}{p}$ to be definably isomorphic to $H(\calC)$, which by Fact \ref{fact: if bind-grp lin then GmGa} gives the corollary ($\bg{F}{p}$ must be of dimension $k$ because it acts regularly on a type of rank $k$).
    
\end{proof}

We are now almost ready to prove the main theorem of this subsection. Before we do so, we recall again the definition of the Lie derivative. The system (\ref{sys of eq}) induces a rational vector field $X$ on $\mathbb{A}^n$, meaning a rational section of the tangent bundle. Given any $g \in F(x_1, \cdots, x_k)$, the Lie derivative measure the variation of $g$ along the flow of that vector field, it is given by:
\[\calL_X(g) = \sum\limits_{i=1}^k \partials{g}{x_i}f_i \text{ .}\]
Using the chain rule, it is straightforward to compute that for any $g \in F(x_1, \cdots , x_k)$ and any solution $(a_1, \cdots ,a_k)$ of (\ref{sys of eq}), we have
\[\delta(g(a_1, \cdots , a_k)) = \calL_X(g)(a_1, \cdots , a_k)\]
\noindent and this will be used in the proof below.

By Lemma \ref{lem: lin-bind-group-im}, if the generic type of an equation of the form (\ref{sys of eq}) is $\calC$-internal and weakly $\calC$-orthogonal, then it has a linear binding group, and Corollary \ref{cor: interdef-full-log} applies. Given this, we can now state and prove the main theorem of this subsection:

\begin{theorem}\label{theo: intern-crit-weak-ortho}
    Let $F$ be an algebraically closed field of constants, some $f_1, \cdots f_k \in F(x_1, \cdots , x_k)$ and $p$ the generic type of the system:
    \[\begin{cases}
    y_1' = f_1(y_1, \cdots , y_k) \\
    \vdots \\
    y_k' = f_k(y_1, \cdots , y_k)
    \end{cases}\]
     \noindent as well as $X$ the associated vector field. Then $p$ is almost $\calC$-internal and weakly $\calC$-orthogonal if and only if there are $g_1 , \cdots ,g_k \in F(x_1, \cdots , x_k)$, algebraically independent over $F$, such that for all $1 \leq i \leq k$, either:
    \begin{itemize}
        \item $\calL_X(g_i) = \lambda_i g_i$ for some $\lambda_i \in F \setminus \{ 0 \}$, 
        \item $\calL_X(g_i) = 1$,
    \end{itemize}
    where the $\lambda_i$ are $\mathbb{Q}$-linearly independent and the second equality is true for at most one $i$.
\end{theorem}

\begin{proof}

    If $p$ is almost $\calC$-internal, it is a well-know fact that there is an $F$-definable map $\pi$ such that $\pi(p)$ is $\calC$-internal and $a \in \acl(\pi(a)A)$ for any $a \models p$. The type $\pi(p)$ is weakly $\calC$-orthogonal because $p$ is, and it has linear binding group by Lemma \ref{lem: lin-bind-group-im}. 
    
    Since $p$ is almost $\calC$-internal, we have $U(p) = k$, and therefore $U(\pi(p)) = k$. Corollary \ref{cor: interdef-full-log} shows that $\pi(p)$ is interdefinable with the generic type of the full logarithmic-differential equation over $F$ on either $G_m^{k-1} \times G_a$ or $(G_m)^k$. Let $\nu : \pi(p) \rightarrow \nu(\pi(p))$ be the $F$-definable map witnessing it.

    Assume first that we are in the $G_m^k$ case. Then by Lemma \ref{lem: fullness for Gm}, the type $\nu(\pi(p))$ is the generic type of a system:
    \[\begin{cases}
        z_1' = \lambda_1 z_1 \\
        \vdots \\
        z_k' = \lambda_k z_k
    \end{cases}\]
    where the $\lambda_j$ are $\mathbb{Q}$-linearly independent elements of $F$. There are $g_1, \cdots ,g_k \in F(x_1, \cdots , x_k)$ such that $\nu \circ \pi$ is given by $\nu(\pi(\abar)) = \left( g_1(\abar), \cdots , g_k(\abar) \right)$ for any $\abar = (a_1 , \cdots , a_k) \models p$ 
    
    For all $1 \leq j \leq k$ we have $\delta(g_j(\abar)) = \lambda_j g_j(\abar)$ by our choice of $g_j$, which gives us:
    \[\lambda_j g_j(\abar) = \delta(g_j(\overline{a}))= \calL_X(g_j)(\abar)\]
    \noindent and since this equation does not involve the derivative, by genericity of $\abar$, we obtain the following equality in $F(x_1, \cdots , x_k)$:
    \[\calL_X(g_j) = \lambda_j g_j \text{ .}\]
    The $G_m^{k-1} \times G_a$ case is similar, the only difference being that we now use Lemma \ref{lem: fullness} and Lemma \ref{lem: fullness for Ga} to see that $\pi(p)$ is interdefinable with the generic type of:
    \[\begin{cases}
        z_1' = \lambda_1 z_1 \\
        \vdots \\
        z_{k-1}' = \lambda_{k-1} z_{k-1} \\
        z_k' = 1
    \end{cases}\]
    \noindent where $\lambda_1, \cdots , \lambda_{k-1} \in F$ are $\mathbb{Q}$-linearly independent. By the same computation as before, we get the partial differential equation for $g_k$. Note that $g_k$ is the only $g_i$ satisfying $\mathcal{L}_X(g_i)=1$.

    The algebraic independence of the $g_i$ over $F$ is a consequence of the Kolchin-Ostrowski theorem (see \cite[top of page 1156]{kolchin1968algebraic}). Indeed, picking some $\abar \models p$, for all $i$ we have that $g_i(\abar)$ is either primitive or exponential over $F$ (in the terminology of \cite{kolchin1968algebraic}). If the $g_i$ were algebraically dependent over $F$, the $g_i(\abar)$ would be as well, and Kolchin-Ostrowski gives a relation of the form
    \[\prod\limits_i g_i(\abar)^{e_i} = \gamma \]
    \noindent for some $e_i \in \mathbb{Z}$ not all zero and some $\gamma \in F \setminus \{ 0 \}$, where the product is taken over all $g_i$ such that $\calL_X(g_i) = \lambda_i g_i$. Since $F$ is a field of constants, taking a derivative yields
    \begin{align*}
        0 & = \left(\prod\limits_i g_i(\abar)^{e_i}\right)' \\
        & = \sum\limits_{j} \left( e_j \lambda_j g_j(\abar)^{e_j} \prod\limits_{i \neq j} g_i(\abar)^{e_i} \right) \\
        & = \gamma \sum\limits_j e_j \lambda_j \text{ .}
    \end{align*}
    \noindent In other words, the $\lambda_i$ are $\mathbb{Q}$-linearly dependent, a contradiction.

    Conversely, suppose that we have some $g_i$ satisfying the statement of the theorem. Assume for now that there is no $i$ such that $\calL_X(g_i) = 1$. Because the $\lambda_i$ are linearly independent over $\mathbb{Q}$, the logarithmic differential equation $\dl_{G_m^k}(x) = (\id_{G_m^k}, (\lambda_1 , \cdots , \lambda_k))$ is full by Lemma \ref{lem: fullness for Gm}, let $q$ be its generic type. 

    We define the map $\pi $ on $p$ by $\pi(\abar) = (g_1(\abar), \cdots , g_k(\abar))$. The same computation as before shows that $\delta(g_j(\abar)) = \lambda_j g_j(\abar)$ for any $1 \leq j \leq k$ and $\abar \models p$. For any $1 \leq j \leq k$, we have that $g_j(\abar) \neq 0$ by genericity of $\abar$. As $q$ is isolated by $\dl_{G_m^k}(x) = (\id_{G_m^k}, (\lambda_1 , \cdots , \lambda_k))$ and $x \in G_m^{k}$, this implies that $\pi(\abar) \models q$ for any $\abar \models p$. 

    Therefore $\pi(p(\calU)) = q(\calU)$, and since the Lascar rank of $q$ is $k$, the Lascar rank of $p$ must be $k$ as well. In particular $\pi$ has finite fibers, and since $q$ is $\calC$-internal, we finally obtain that $p$ is almost $\calC$-internal. By fullness $q$ is weakly $\calC$-orthogonal, and therefore $p$ also is weakly $\calC$-orthogonal as $p$ and $q$ are interalgebraic.

    If one of the $g_i$ satisfies $\calL_X(g_i) = 1$, the same proof works \emph{mutatis mutandis}.
     
\end{proof}

\begin{remark}
    \sloppy It is natural to ask if we could characterize outright $\calC$-internality. It should be equivalent to the map $(x_1, \cdots, x_k) \rightarrow (g_1(x_1, \cdots , x_k), \cdots , g_k(x_1, \cdots , x_k))$ inducing a birational transformation of $\mathbb{A}^k$.
\end{remark}

The same ideas easily give a criteria for weak orthogonality to the constants, which is implicit in the literature but, as far as the authors know, was not stated previously in this form. Note that we are not necessarily working over constant parameters here, and thus this is not an immediate corollary of the theorem. Recall that if $F$ is any algebraically closed differential field, there is a unique derivation $\delta$ on $F(x)$ extending that of $F$ and with $\delta(x) = 0$. If $g \in F(x)$, we denote $\delta(g)= g^{\delta}$.

\begin{proposition}\label{prop: weak-ortho-crit}
   Let $F$ be an algebraically closed differential field, some rational functions $f_1, \cdots , f_k \in F(x_1, \cdots , x_k)$ and $p$ the generic type of the system:
    \[\begin{cases}
    y_1' = f_1(y_1, \cdots , y_k) \\
    \vdots \\
    y_k' = f_k(y_1, \cdots , y_k)
    \end{cases}\]
    Then $p$ is not weakly orthogonal to $\calC$ if and only if there is a non-constant $g \in F(x_1, \cdots , x_k)$ such that $\sum\limits_{i=1}^{k} \partials{g}{x_i}f_i  + g^{\delta} = 0$.
\end{proposition}

\begin{proof}
    First suppose that such a $g$ exists. Consider some $\abar = (a_1, \cdots , a_k) \models p$, we compute that:
    \[\delta(g(\abar)) = \sum\limits_{i=1}^{k} \partials{g}{x_i}(\abar)f_i(\abar) + g^{\delta}(\abar) = 0\]
    \noindent so $g(\abar) \in \calC$. Moreover $g(\abar) \not\in F$ by genericity of $\abar \models p$ (and because $g$ is not constant). Thus $g(\abar)$ and $\abar$ fork over $F$, which implies that $p$ is not weakly orthogonal to $\calC$.

    For the converse, assume that $p$ is not weakly orthogonal to $\calC$. By \cite[Lemma 2.1]{freitag2023bounding}, there is $\abar = (a_1, \cdots , a_k) \models p$ and $b \in \calC$ with $b \in \dcl(\abar F) \setminus F$. This implies that there is $g \in F(x_1, \cdots , x_k)$ such that $b = g(\abar)$. By the same computation as before and genericity of $\abar$, we obtain the desired equation.
\end{proof}

We remark that Theorem \ref{theo: intern-crit-weak-ortho} also gives an easy criteria for orthogonality to the constants:

\begin{corollary}\label{cor: ortho crit}
   Let $F$ be an algebraically closed field of constants, some $f_1, \cdots f_k \in F(x_1, \cdots , x_k)$ and $p$ the generic type of the system:    
    \[\begin{cases}
    y_1' = f_1(y_1, \cdots , y_k) \\
    \vdots \\
    y_k' = f_k(y_1, \cdots , y_k)
    \end{cases}\]
    \noindent as well as $X$ the associated vector field. Then $p$ is not orthogonal to the constants if and only if there is $g \in F(x_1, \cdots , x_k) \setminus F$ and some $\gamma \in \calC$ such that either $\calL_X(g) = \gamma g$ or $\calL_X(g)  = \gamma$. 
\end{corollary}

\begin{proof}
    If such a $g$ exists, then it is easy to see that for any $\abar = (a_1, \cdots , a_k) \models p$, the type $\tp(g(\abar)/F)$ is $\calC$-internal, because it satisfies a logarithmic differential equation on either $G_a$ or $G_m$. Therefore $g(p)$ is not orthogonal to $\calC$ and neither is $p$.

    Conversely, suppose that $p$ is not orthogonal to $\calC$. Then by \cite[Chapter 7, Corollary 4.6]{pillay1996geometric}, there is some $\abar = (a_1, \cdots , a_k) \models p$ and some $b \in \dcl(\abar F) \setminus F$ such that $q = \tp(b/F)$ is $\calC$-internal. Thus we have an $F$-definable function $h : p \rightarrow q$. If $q$ is not weakly $\calC$-orthogonal, then neither is $p$, and we conclude by Proposition \ref{prop: weak-ortho-crit}. So we may assume that $q$ is weakly $\calC$-orthogonal. By Lemma \ref{lem: lin-bind-group-im}, the binding group of $q$ is linear. By the proof of Theorem \ref{theo: intern-crit-weak-ortho}, we also have an $F$-definable function from $q$ to the generic type of a full logarithmic differential equation over either $G_m^l$ or $G_m^{l-1} \times G_a$ for some $l \in \mathbb{N}$. Composing these two functions and taking a coordinate function, we get the result using the same argument as in the proof of Theorem \ref{theo: intern-crit-weak-ortho}.
\end{proof}

\subsection{Non-weakly orthogonal case}\label{subsec: not weak orth}

In this section, we prove a version of Theorem \ref{theo: intern-crit-weak-ortho}, without assuming that the type is weakly $\calC$-orthogonal. We give a separate proof as a number of new issues appear. This involves some repetition, which we tried to keep to a minimum.

The key point is that non weak orthogonality to $\calC$ is always witnessed by a function with weakly $\calC$-orthogonal fibers:

\begin{lemma}\label{lem: weak-ortho-cored}
    Let $F$ be an algebraically closed differential field of constants and $p \in S(F)$ a finite rank type. Then there is an $F$-definable map $\pi : p \rightarrow \pi(p)$ such that $\pi(p)(\calU) \subset \calC$ and for any $a \models p$, the type $\tp(a/\pi(a)F)$ is stationary and weakly $\calC$-orthogonal. 
\end{lemma}

\begin{proof}
    As $\calC$ is stably embedded, we have, for any $a \models p$, that $\tp(a/ \dcl(aF) \cap \calC)$ isolates $\tp(a/\calC )$ (see the appendix of \cite{chatzidakis1999model}). The set $\dcl(aF)$ is a finitely generated field extension of $F$ because $\tp(a/F)$ has finite rank. The set $\dcl(aF) \cap \calC$ is a definably closed subset of the stably embedded pure algebraically closed field $\calC$, and therefore is a field. In particular it is a subextension of the finitely generated extension $F < \dcl(aF)$, which implies that it is also finitely generated. 
    
    So there are finitely many $F$-definable functions $\pi_1, \cdots , \pi_n$, all defined on $p$ and with image in $\calC$, such that $\dcl(aF) \cap \calC = \dcl(\pi_1(a), \cdots , \pi_n(a),F)$. Letting $\pi = (\pi_1, \cdots, \pi_n)$, we obtain that $\tp(a/\pi(a)F) \models \tp(a/\calC)$. 

    Note that $\tp(a/\calC)$, as $\calC$ is algebraically closed, must be stationary. Therefore $\tp(a/\pi(a)F)$ also is stationary. Moreover, if $d$ is any tuple of constants, then $\tp(a/\pi(a)F) \models \tp(a/\pi(a)Fd)$ and in particular $a \forkindep_{\pi(a)F} d$. So $\tp(a/\pi(a)F)$ is weakly orthogonal to $\calC$. 
\end{proof}

We can now prove our main theorem:
\begin{theorem}\label{theo: inter-crit-not-weak-ortho}
    Let $F$ be an algebraically closed field of constants, some $f_1, \cdots f_k \in F(x_1, \cdots , x_k)$ and $p$ the generic type of the system:
    \[\begin{cases}
    y_1' = f_1(y_1, \cdots , y_k) \\
    \vdots \\
    y_k' = f_k(y_1, \cdots , y_k)
    \end{cases}\]
    \noindent as well as $X$ the associated vector field. Then $p$ is almost $\calC$-internal if and only if there are $g_1 , \cdots ,g_k \in F(x_1, \cdots , x_k)$, algebraically independent over $F$, such that for all $1 \leq i \leq k$, either:
    \begin{itemize}
        \item $\calL_X(g_i) = \lambda_i g_i$ for some $\lambda_i \in F$,
        \item $\calL_X(g_i) = 1$.
    \end{itemize}

    Moreover, the type $p$ is weakly $\calC$-orthogonal if and only if the $\lambda_i$ are $\mathbb{Q}$-linearly independent and at most one $g_i$ satisfies $\mathcal{L}_X(g_i) = 1$.
\end{theorem}

\begin{proof}
    \sloppy We start with the left to right direction, and assume that $p$ is almost $\mathcal{C}$-internal. There is a finite-to-one $F$-definable map $\rho:p\to q$ such that $q$ is $\mathcal{C}$-internal.

    Let $\pi$ be the map from Lemma \ref{lem: weak-ortho-cored} applied to $q$. For any $\abar \models q$, we denote $F_{\pi(\abar)} : = F(\pi(\abar))$ the field generated by $\pi(\abar)$ over $F$ and $F_{\pi(\abar)}^{\mathrm{alg}} = F(\pi(\abar))^{\mathrm{alg}}$ its algebraic closure. For any $\abar \models q$, the type $q_{\pi(\abar)} := \tp(\abar/ F_{\pi(\abar)})$ is stationary, weakly $\calC$-orthogonal and $\calC$-internal. Since we will want to work over algebraically closed sets of parameters, we will also consider its unique extension $\overline{q_{\pi(\abar)}} := \tp(\abar/ F_{\pi(\abar)}^{\mathrm{alg}})$.

    By \cite[Lemma 2.7]{jaoui2022abelian}, the binding group of $q$ is definably isomorphic to the binding group of $\tp(\abar_1, \cdots, \abar_m/F\pi(\abar_1), \cdots , \pi(\abar_m))$, where $\abar_1, \cdots , \abar_m$ is a fundamental system of solutions of $q$. This binding group is a definable subgroup of the cartesian product $\prod\limits_{i=1}^m \bg{F_{\pi(\abar_i)}}{q_{\pi(\abar_i)}}$. By Lemma \ref{lem: lin-bind-group-im}, the binding group of $q$ is definably isomorphic to the $\calC$-points of a linear algebraic group, and therefore it is also the case for each $\bg{F_{\pi(\abar_i)}}{q_{\pi(\abar_i)}}$.

    Therefore the binding group $\bg{F_{\pi(\abar)}}{q_{\pi(\abar)}}$ is definably isomorphic to the $\calC$-points of a linear algebraic group for all $\abar \models q$. Fix such a realisation $\abar$. By \cite[Lemma 2.1]{jaoui2022abelian}, the binding group $\bg{F_{\pi(\abar)}^{\mathrm{alg}}}{\overline{q_{\pi(\abar)}}}$ is a definable normal subgroup of $\bg{F_{\pi(\abar)}}{q_{\pi(\abar)}}$ and thus is also definably isomorphic to the $\calC$-points of a linear algebraic group. 

    \sloppy By Corollary \ref{cor: interdef-full-log}, the type $\overline{q_{\pi(\abar)}}$ is fundamental and interdefinable with the generic type of a full logarithmic differential equation over $F_{\pi(\abar)}^{\mathrm{alg}}$ on either $G_a \times (G_m)^{l-1}$ or $(G_m)^l$, where $l = U(\overline{q_{\pi(\abar)}})$. The map witnessing interdefinability is given by $(\gamma_1(y_1, \cdots , y_k), \cdots , \gamma_l(y_1,\cdots , y_k))$, where $\gamma_i$ is an $F_{\pi(\abar)}^{\mathrm{alg}}$-definable map, and each $\gamma_i$ maps $\overline{q_{\pi(\abar)}}$ to solutions of a full logarithmic differential equation on $G_m$ or $G_a$ over $F_{\pi(\abar)}^{\mathrm{alg}}$.

    For each $i$, fix $\lambda_i$ such that $\gamma_i$ maps $\overline{q_{\pi(\abar)}}$ to solutions of either $z' = \lambda_i z$, or $z ' = \lambda_i$ if $i = l$. Note that none of the $\lambda_i$ are zero, as $\overline{q_{\pi(\abar)}}$ is weakly orthogonal to $\mathcal{C}$. In fact, they are $\mathbb{Q}$-linearly independent by Lemma \ref{lem: fullness}. As discussed in Lemma \ref{lem: fullness for Ga}, if $\gamma_l$ maps to solutions of $z ' = \lambda_l$, we can compose with the map $z \rightarrow \frac{z}{\lambda_l}$ to obtain solutions of $z' = 1$ instead. We make this harmless assumption for the rest of the proof.
    
    We now show that in the $(G_m)^l$ case, we can obtain that $\lambda_i \in F$:

    \begin{claim}
        For all $i$, we have $\lambda_i \in F$
    \end{claim}

    \begin{proof}[Proof of claim]
        By the previous discussion, we can assume that $\gamma_i$ maps $\overline{q_{\pi(\abar)}}$ to solutions of $z' = \lambda_i z$. Let $d = \gamma_i(\overline{a})$. Since $d$ is algebraic over $F( \overline{a})$, the type $\tp(d/F)$ is almost $\mathcal{C}$-internal. 

        Note that $U(\tp(d/F(\lambda_i))) = 1$, which implies that $\tp(d/F(\lambda_i)) = \dl^{-1}(\tp(\lambda_i/F))$, which is defined as the type of some (any) $w$ such that $\dl(w) \models \tp(\lambda_i/F)$ and $w \not \in \acl(\dl(w)F)$ (see \cite[Page 5]{jin2020internality} for more details).

        Let $r$ be the generic type of $\calC$ over $\mathbb{Q}^{\mathrm{alg}}$, it is well-known (see \cite[Lemma 4.2]{chatzidakis2015differential} for a proof) that $\dl^{-1}(r)$ is not almost $\calC$-internal. Assuming that $\lambda_i \not\in F$, we obtain $\tp(\lambda_i/F) = r\vert_F$ because $F$ is algebraically closed, and therefore $\tp(d/F) = \dl^{-1}(r)\vert_F$ is not almost $\calC$-internal either. But this is a contradiction, as $\tp(d/F)$ is almost $\mathcal{C}$-internal.
    \end{proof}

    Using this, we can make the functions $\gamma_i$ descend to $F$:

    \begin{claim}
        There are integers $n_i$ and $F$-definable maps $\tilde{\gamma_i}$, for $i \in \{1, \cdots l \}$ mapping $q$ to the generic type of $z' = n_i \lambda_i z$ if $i < l$, and either $z' = n_l \lambda_l z$ or $z ' = 1$.
    \end{claim}

    \begin{proof}[Proof of claim]

        Fix some $i$, and assume first that $\gamma_i$ maps $\overline{q_{\pi(\overline{a})}}$ to solutions of a logarithmic differential equation on $G_m$. Since $q_{\pi(\abar)} \models \overline{q_{\pi(\abar)}}$, the function $\gamma_i$ defines a map on $q_{\pi(\abar)}$.
        
        The $F_{\pi(\overline{a})}^{\mathrm{alg}}$-definable function $\gamma_i$ has finitely many images $\gamma_{i,1}, \cdots , \gamma_{i,n_i}$ under $\Aut_{F_{\pi(\abar)}}(\calU)$, all mapping $q_{\pi(\abar)}$ to the generic type of $z ' = \lambda_i z$, since $\lambda_i \in F$ by the previous claim.

        Consider $\hat{\gamma_i} = \gamma_{i,1} \times \cdots \times \gamma_{i,n_i}$. We check that $\hat{\gamma_i}(\overline{a})$ is not a constant. Otherwise, it would be a constant solution of $z' = n_i \lambda_i z$. Since $\lambda_i \neq 0$ and $n_i \neq 0$, this would imply that $\hat{\gamma_i}(\overline{a}) = 0$, which implies that $\gamma_{i,j}(\overline{a}) = 0$ for some $1 \leq j \leq n_i$. But this is impossible as $\gamma_{i,j}(\overline{a})$ realizes the generic type of $z' = \lambda_i z$.

        From the previous discussion, we obtain that $\hat{\gamma_i}(\overline{a})$ realizes the generic type of $z' = n_i \lambda_i z$ over $F$, as all non-generic solutions must be constant. Thus $\hat{\gamma_i}$ maps $q_{\pi(\abar)}$ to the generic type of $z' = n_i \lambda_i z$. The new map $\hat{\gamma_i}$ is fixed by $\Aut_{F_{\pi(\abar)}}(\calU)$ and therefore $F_{\pi(\abar)}$-definable. 

        In the other case, we have $i = l$ and $\gamma_l$ maps to the generic type of $z' = 1$. We can consider its conjugates $\gamma_{l,1}, \cdots , \gamma_{l,n_l}$ under $\Aut_{F_{\pi(\abar)}}(\calU)$ and replace it with the map $\hat{\gamma_l} = \gamma_{l,1} + \cdots + \gamma_{l,n_l}$, which an is $F_{\pi(\overline{a})}$-definable map from $q_{\pi(\abar)}$ to the generic type of $z' = n_l$.  By Lemma \ref{lem: fullness for Ga} again, we can replace $n_l$ by $1$. 

        So $\hat{\gamma_i}$ is $ F_{\pi(\abar)}$-definable for all $i$. Moreover, we know that $\pi(\abar) \in F(\abar)$, and therefore $\pi(\overline{a})$ is defined by some $F$-formula $\theta$. We can replace $\pi(\overline{a})$ in the formulas defining the $\hat{\gamma_i}$ by $\theta(\overline{a})$ to obtain some $F$-definable map $\tilde{\gamma_i}$, defined on $q$. 
        
    \end{proof}

    Recall that the $\lambda_i $ are $\mathbb{Q}$-linearly independent. For all $i$, let $\tilde{\lambda}_i=n_i\lambda_i$,  the $\tilde{\lambda_i}$ are also $\mathbb{Q}$-linearly independent. To summarize, have obtained $F$-definable maps $\tilde{\gamma_1}, \cdots, \tilde{\gamma_{l-1}} , \tilde{\gamma_l}$ such that $\tilde{\gamma}= (\tilde{\gamma_1}, \cdots , \tilde{\gamma_l})$ maps $q$ to the generic type of:
    \[\begin{cases}
        z_1' = \tilde{\lambda}_1 z_1 \\
        \vdots \\
        z_{l-1}' = \tilde{\lambda}_{l-1} z_{l-1}\\
        z_l' = \tilde{\lambda}_l z_l \text{ or } z_l' =1
    \end{cases}\]
    where the $\tilde{\lambda}_i$ are $\mathbb{Q}$-linearly independent elements of $F$, and the last line depends on whether the binding group of $\overline{q_{\pi(\abar)}}$ is $G_m^l$ or $G_m^{l-1} \times G_a$. Using a similar proof to that of Theorem \ref{theo: intern-crit-weak-ortho} replacing the $\pi:p\to\pi(p)$ and $\nu:\pi(p)\to\nu(\pi(p))$ in Theorem \ref{theo: intern-crit-weak-ortho} with $\rho:p\to q$ and $\tilde{\gamma} : q\to \tilde{\gamma}(q)$ respectively, we obtain rational functions $g_1, \cdots, g_l$, algebraically independent over $F$, such that for all $i$, we have $\calL_X(g_i) = \lambda_i g_i$ or $= 1$. 

    To obtain the rest of the $g_i$, consider the map $\pi \circ \rho : p \rightarrow \pi \circ \rho (p)$. Let $\overline{b} \models p$ and $\overline{a} = \rho(\overline{b})$. Recall that $\pi(\overline{a})$ is a tuple of constants, and thus by quantifier elimination, there are $h_1, \cdots , h_n \in F(x_1, \cdots , x_k)$ such that $\pi\circ \rho(\bbar) = (h_1(\bbar), \cdots , h_n(\bbar))$. Since $p$ is almost $\mathcal{C}$-internal, we know that $U(p)$ is equal to the transcendence degree of any realization over $F$, so $U(p) = k$. We can compute that:
    \begin{align*}
        l & = U(\overline{b}/\pi \circ \rho(b)F) \\
        & = U(\overline{b}/F) - U(\pi\circ \rho(\overline{b})/F) \\
        & = k - U(\pi\circ \rho(\overline{b})/F)
    \end{align*}
    \noindent Therefore a maximal algebraically independent tuple from $\pi \circ \rho(\overline{b})$ has length $k-l$, and we may assume, without loss of generality, that $h_1(\overline{b}), \cdots, h_{k-l}(\overline{b})$ are algebraically independent over $F$. In particular the rational functions $h_1, \cdots , h_{k-l}$ are algebraically independent over $F$. Note that since $\tilde{\gamma} \circ \rho (p) = \tilde{\gamma}(q)$ is weakly $\mathcal{C}$-orthogonal, we have that $\tilde{\gamma} \circ \rho (\overline{b})$ and $\pi \circ \rho(\overline{b})$ are independent over $F$. This implies that the rational functions $g_1, \cdots , g_{k-l}, h_1, \cdots , h_l$ are algebraically independent over $F$. We pick the rest of the $g_i$ to be $h_1, \cdots , h_{l-k}$, it is again easy to get $\calL_X(g_i) = 0$ for all $i > l$.  

    Conversely suppose that there are $g_1, \cdots , g_{k}$ satisfying the assumption of the theorem. Assume for now that all $g_i$ satisfy $\calL_X(g_i) = \lambda_i g_i$ for some $\lambda_1, \cdots , \lambda_{k}$. 
    
    By reordering, we may assume that $\lambda_1, \cdots , \lambda_l$ is a maximal $\mathbb{Q}$-linearly independent subset of $\{ \lambda_1 ,\cdots , \lambda_k\}$. The logarithmic differential equation $\dl_G(x) = (\id_G,(\lambda_1, \cdots , \lambda_{l}))$ is full. Let $q$ be the generic type of its solutions. Using the same proof as in Theorem \ref{theo: intern-crit-weak-ortho}, we obtain that $g : \abar \rightarrow (g_1(\abar), \cdots , g_{l}(\abar))$ defines a surjective map from $p(\calU)$ to $q(\calU)$. 

    Pick $\lambda_j \in \{ \lambda_{l+1}, \cdots , \lambda_k\}$. There are $a_{j,1}, b_{j,1}, \cdots , a_{j,l},b_{j,l} \in \mathbb{Z}$ such that $\lambda_j = \sum\limits_{i = 1}^l \frac{a_{j,i}}{b_{j,i}} \lambda_i$. From this we deduce that $h_j = \frac{\prod\limits_{i=1}^l g_i^{\left( \prod\limits_{s \neq i} b_{j,s}\right)a_{j,i}}}{g_j^{\prod\limits_{i=1}^l b_{j,i}}}$ maps $p$ to the constants. 
    
    Since $g_j$ and $h_j$ are interalgebraic over $F(g_1, \cdots , g_l)$ for all $j \geq l+1$, we get that $g_1, \cdots , g_l, h_{l+1}, \cdots h_k$ are algebraically independent over $F$. Now define $\pi : \overline{a} \rightarrow (h_{l+1}(\overline{a}), \cdots , h_k(\overline{a}))$. Because the $h_j$ are algebraically independent over $F$, this defines a map from $p$ to $r^{(k-l)}$, where $r$ is the type of a constant transcendental over $F$. Note that since the logarithmic differential equation $\dl_G(x) = (\id_G,(\lambda_1, \cdots , \lambda_{l}))$ is full, the type $q$ is weakly $\calC$-orthogonal. Therefore, for any $\abar \models p$, we have that $U(\pi(\abar), g(\abar)/F)  = U(\pi(\abar)/F) + U(g(\abar)/F)  = k$.

    Again we remark that the $U$-rank of $p$ must be less or equal to $k$, and because $g \times \pi$ defines a surjective map to the set of realizations of a type of $U$ rank $k$, we must have $U(p) = k$. In particular the fibers of $g \times \pi $ must be finite. As $g(p) \otimes \pi(p)$ is $\calC$-internal, the type $p$ is almost $\calC$-internal.

    Now assume that some $g_i$ satisfies $\mathcal{L}_X(g_i) = 1$. Note that if $g_j$ also satisfies $\mathcal{L}_X(g_j) = 1$, then $\mathcal{L}_X(g_i-g_j) = 0$. We can thus again, by reordering, pick a maximal set $\{ \lambda_1 , \cdots , \lambda_l \}$ of $\mathbb{Q}$-linearly independent $\lambda_i$, then pick $g_{l+1}$ with $\mathcal{L}_X(g_{l+1}) = 1$. The rest of the $g_j$ can combined to obtain some $h_j$ mapping $p$ to the constants, and we prove that $g_1, \cdots , g_{l+1}, h_{l+1}, \cdots h_k$ are algebraically independent as was done in the previous paragraph. We can then conclude as before.  

    The left to right implication of the moreover part follows from Theorem \ref{theo: intern-crit-weak-ortho}, the right to left is a consequence of the above proof: if the $\lambda_i$ are $\mathbb{Q}$-linearly independent and at most one $g_i$ satisfies $\mathcal{L}_X(g_i) = 1$, then $g$ is a finite-to-one map from $p$ to the generic type of a full logarithmic differential equation, and thus must be weakly $\mathcal{C}$-orthogonal.
     
\end{proof}

\section{Applications}\label{sec: applications}

From Theorem \ref{theo: intern-crit-weak-ortho} and Theorem \ref{theo: inter-crit-not-weak-ortho}, we see that to understand whether or not the generic type of the solution set of an autonomous algebraic ordinary differential equations is almost $\calC$-internal, we have to understand all possible maps to solutions of logarithmic differential equations on $G_m$ or $G_a$. In turn, this is done by determining rational function solutions to certain order one partial differential equations. In this section, we apply these theorems to specific differential equations from the literature. 

We will make heavy use of Laurent series methods, and we thus recall a few useful facts here. Let $K$ be any field. Then the field of rational functions $K(x)$ embeds naturally in the field of formal Laurent series $K((x))$, which consists of all formal sums of the form $\sum\limits_{i= N}^{\infty} a_i x^i$, with $N \in \mathbb{Z}$, all $a_i$ elements of $K$ and $a_N \neq 0$. Moreover, if $K(x)$ is equipped with a derivation $\delta$, then it extends naturally to $K((x))$. We will often apply this to $K(x_1,x_2)$ equipped with the derivations $\partials{}{x_1}$ and $\partials{}{x_2}$, which we will view as a differential subfield of either $K(x_1)((x_2))$ or $K(x_2)((x_1))$.

\subsection{Poizat equations}\label{subsec: poizat}

In \cite{freitag2023equations}, Freitag, Jaoui, Marker and Nagloo consider Poizat equations:
\[y'' = y'f(y)\]
\noindent where $f \in F(x)$ and $F$ is a field of constants. They show that its set of solutions is strongly minimal if and only if $f$ is not the derivative of some $g \in F(x)$. Otherwise, they show that there is a definable map $\pi : y \rightarrow y'-g(y)$ from realizations of $p$ to the constants. In particular $\pi(p)$ is internal to the constants. Given this, there are three cases for the semi-minimal analysis of $p$: 

\begin{enumerate}
    \item $p$ is almost internal to the constants, for example if $g = c \in F$,
    \item the fibers of $\pi$ are almost internal to the constants, but $p$ is not almost internal to the constants, for example if $g(y) = y$, by \cite[Lemma 7.8]{freitag2023equations},
    \item the fibers of $\pi$ are orthogonal to the constants. By \cite[Corollary 7.12]{freitag2023equations}, this is the case if there is a degree three polynomial $P \in F[x]$ such that $\partials{P}{x} = f$.   
\end{enumerate}

The authors ask, in \cite[Question 7.6]{freitag2023equations}, whether $y'' = cy'$ is the only case for which $p$ is almost internal to the constants. We show that this is indeed the case. Fix some algebraically closed field of constants $F$.

\begin{proposition}
    Let $p$ be the generic type of  $y''=y'f(y)$ for some $f \in F(x)$, which is interdefinable with the generic type of the system
    \[\begin{cases}
    y_1' = y_2 \\
    y_2' = y_2f(y_1) \text{ .}
    \end{cases}\]
    Then $p$ is almost $\calC$-internal if and only if $f \in F$. 
\end{proposition}

\begin{proof}
    Suppose that $p$ is almost $\calC$-internal. By the previous discussion, there is some $g \in F(x)$ such that $\partials{g}{x} = f$. In particular, as was observed in \cite{freitag2023equations}, we obtain an $F$-definable map $ \pi : y \rightarrow y' - g(y)$ from $p(\calC)$ to the constants. 

    By Theorem \ref{theo: inter-crit-not-weak-ortho}, to determine whether $p$ is almost $\calC$-internal, we must determine if there can be two $F$-definable maps to either the constants, or solutions of logarithmic differential equations on $G_a$ or $G_m$. We already have a map $\pi$ to the constants, and as $p$ is the generic type of an order two differential equation, we must have $p(\calU) \not\subset \calC$. In particular, if $p$ is almost $\calC$-internal, the second map cannot go to the constants. Therefore we only have to consider the $G_a$ and $G_m$ cases. 

    In the $G_m$ case, we must solve: 

    \[ \partials{h}{x_1}x_2+ \partials{h}{x_2}x_2f(x_1) = \lambda h\]

    \noindent for $\lambda \in F \setminus \{ 0\}$ and $h \in F(x_1,x_2)$. We view $h$ as a Laurent series $h = \sum\limits_{i=N}^{\infty} a_i x_2^i$ with $a_i \in F(x_1)$ and $N$ being least such that $a_N \neq 0$. Looking at coefficients for $x_2^N$, we obtain the equation $Na_Nf(x_1) = \lambda a_N$ (as $a_{N-1} = 0$). If $N \neq 0$, since $a_N \neq 0$, we get $f(x_1) = \frac{\lambda}{N}$, i.e. $f \in F$. If $N = 0$, the equation becomes $\lambda a_N = 0$, which is impossible as both $a_N$ and $\lambda$ are non-zero. 

    In the $G_a$ case, we must solve:

    \[ \partials{h}{x_1}x_2+ \partials{h}{x_2}x_2f(x_1) = 1 \]

    \noindent for $h \in F(x_1,x_2)$. We again view $h$ as a Laurent series. If $N \neq 0$, looking at the $x_2^N$ coefficients gives $Na_N f(x_1) = 0$, therefore $f(x_1) = 0$, so $f \in F$. If $N = 0$, then we get $0 = 1$, a contradiction.

\end{proof}

\subsection{Lotka-Volterra equations}\label{subsec: LV}

A Lotka-Volterra system is given by:

\[
\begin{cases}
    x' = ax + b x y \\
    y' = cy + d x y
\end{cases}
\]
\noindent where $a,b,c,d $ are some constants. We examine this system in light of Theorem \ref{theo: inter-crit-not-weak-ortho}. 

These equations were first suggested by Lotka \cite{lotka1910contribution} (in a slightly different form) for the study of autocatalytic chemical reactions, and later independently by Volterra \cite{volterra1928variations} to model the dynamics of predator-prey biological systems. 

In a biological context, it is further assumed that $a,b,c,d$ are real numbers. Moreover, if $x$ represents the population of prey species, then $a$ is positive because its population, if left alone, grows, and $b$ is negative, because the presence of predators tends to make the population of the prey species decrease. If $y$ is the population of the predator species, similar biological reasoning imposes $c < 0$ and $d>0$ (the predator species needs to eat the prey to survive). 

Given their importance in modeling natural phenomena, finding exact solutions to Lotka-Volterra systems, and variations of them, has been an important problem. As illustrations, see \cite{cherniha2022construction} for a recent survey of the known exact solutions of autonomous Lotka-Volterra systems, and \cite{petrovskii2005exact} for recent work on finding exact solutions to non-autonomous Lotka-Volterra systems. Our works here rules out a large class of potential exact solutions, unless $a = c$, in which case exact solutions were found by Varma \cite{varma1977exact}. This is explained in Corollary \ref{cor: LV-not-classical}.

We will not make any of these assumptions here, but we will assume that none of the parameters $a,b,c$ and $d$ are zero. Using the change of variables $u = \frac{d}{c}x $ and $v = \frac{b}{c} y$, we see that the solution set of this system is interdefinable with the solution set of:
\[
\begin{cases}
    u' = au + cuv \\
    v' = cv + cuv
\end{cases}
\]
\noindent and we can therefore work with that second system instead. We let $p$ be its generic type. 

Let $F = \mathbb{Q}(a,c)^{\mathrm{alg}}$ and $\mu = \frac{a}{c}$. By Theorem \ref{theo: inter-crit-not-weak-ortho}, we need to understand solutions $g \in F(x_0,x_1)$ of the partial differential equations:
\[
\calL_X(g) = c\partials{g}{x_0}\left(\mu x_0 + x_0 x_1\right) + c \partials{g}{x_1} \left( x_1 + x_0 x_1 \right) = 
\begin{cases}
    0 \\
    1 \\
    \lambda g
\end{cases}
\]
\noindent for $\lambda \in F \setminus \{ 0 \}$, where $X$ is the associated vector field. 

We will make use of the following well-known fact throughout the rest of this subsection: a rational function is of the form $\frac{h'}{h}$, for some $h\in F(x)$, if and only if it has no zeros and only poles of order one with integer residues. 

We first deal with the top equation, showing that $p$ is weakly $\calC$-orthogonal. Note that this is equivalent to showing that the Lotka-Volterra system has no rational first integral. This is to be contrasted with the classical fact that $dx+c\ln(x)-by-a\ln(y)$ is a first integral. 

We can divide by $c$ on both sides so we just have to prove:

\begin{lemma}\label{lem: LV weak ortho}
    There is no non-constant solution $g \in F(x_0,x_1)$ to $\partials{g}{x_0}(\mu x_0 + x_0 x_1) + \partials{g}{x_1}(x_1 + x_0 x_1) = 0$. In particular $p$ is weakly $\calC$-orthogonal.
\end{lemma}

\begin{proof}
    Let $g$ be a solution, we show that $g$ must be constant. Write $g$ as an element of $F(x_1)((x_0))$, i.e. there are $N \in \mathbb{Z}$ and some $a_i \in F(x_1)$ such that $g = \sum\limits_{i=N}^{\infty} a_i x_0^i$. 

    Inserting this into the partial differential equation, we obtain, by looking at the order $N$ coefficient:
    \[N(\mu+x_1)a_N + x_1 a_N' = 0.\]
    \noindent When $N\neq 0$ this can be rewritten, since $a_N \neq 0$, as:
    \[ \frac{a_N'}{a_N} = \frac{-N\mu}{x_1} - N\]
    \noindent which does not have any solution in $F(x_1)$.

    Therefore $N=0$ and the equation becomes $x_1a_0'=0$, which implies $a_0 ' =0$, or equivalently $a_0 \in F$. For any $i \geq 1$ we have:
     \[i(\mu + x_1)a_i + x_1 a_i '+ x_1 a_{i-1}' = 0\]
    \noindent In particular, suppose that $a_{i-1}' = 0$, then if $a_i \neq 0$ we have:
    \[ \frac{a_i '}{a_i} = -\frac{i \mu}{x_1} - 1\]
    \noindent which has no solution in $F(x_1)$. Thus if $a_{i-1}' =0$ then $a_i = 0$ for all $i \geq 1$, and we conclude that $g = a_0 \in F$.

\end{proof}

Eliminating the second possibility is similar. Again we divide by $c$ on both sides.

\begin{lemma}\label{lem: LV no Ga}
    There is no solution $g \in F(x_0,x_1)$ to $\partials{g}{x_0}(\mu x_0 + x_0 x_1) + \partials{g}{x_1}(x_1 + x_0 x_1) = \frac{1}{c}$.
\end{lemma}

\begin{proof}
    Let $g$ be a solution. We write $g$ as an element of $F(x_1)((x_0))$, i.e. there $N \in \mathbb{Z}$ and some $a_i \in F(x_1)$ such that $g = \sum\limits_{i=N}^{\infty} a_i x_0^i$, with $a_N \neq 0$.

    Inserting this into the partial differential equation, we obtain: 
\[\begin{cases}
    N(\mu+x_1)a_N+x_1a'_N=0 & \text{ if }N\neq 0\\
    x_1a'_N=\frac{1}{c} & \text{ if } N=0\\
\end{cases}\]
When $N\neq 0,$ we have $\frac{a'_N}{a_N}=\frac{-N\mu}{x_1}-N$ which again has no  solutions in $F(x_1).$

If $N=0$ we obtain $a_0'=\frac{1}{c x_1}$ which has no solution in $F(x_1)$ as it has a non-zero residue.

\end{proof}

We now deal with the $\lambda g$ equation, which is more involved. As before, we divide by $c$ on both sides. Let $\lambda \in F \setminus \{ 0 \}$ and $g$ be a non-constant solution of $\partials{g}{x_0}(\mu x_0 + x_0 x_1) + \partials{g}{x_1}(x_1 + x_0 x_1) = \lambda g$ (where are looking for any $\lambda \in F$, so we can replace $\frac{\lambda}{c}$ by $\lambda$).

We can write $g$ as an element of $F(x_1)((x_0))$, i.e. there are $N \in \mathbb{Z}$ and $a_i \in F(x_1)$ such that $g = \sum\limits_{i= N}^{\infty} a_i x_0^i$ and $a_N \neq 0$. Inserting this into the partial differential equation, we obtain, for all $i$:
\[i(\mu + x_1)a_i + x_1 a_i'+x_1a_{i-1}' = \lambda a_i\]
\noindent where the derivatives are taken with respect to $x_1$. For $i = N$, we obtain:
\[N(\mu + x_1)a_N + x_1a_N' = \lambda a_N\]
\noindent which can be rewritten as:
\[\frac{a_N '}{a_N} = \frac{\lambda -N \mu}{x_1} - N \text{ .}\]
\noindent
Since a rational function can be a logarithmic derivative only if it has no zeroes and only simple poles with integer residues, this implies that $N = 0$ and $\lambda \in \mathbb{Z}$. Note that if $g$ is a solution for $\lambda$, then $g^k$ is a solution for $k \lambda$, for any $k \in \mathbb{Z}$. Therefore, replacing $g$ by one of its powers, we may assume that $\lambda$ is a large positive integer, and we do so for the rest of the proof.

So $a_N = a_0$ is a solution of $\frac{a_0'}{a_0} = \frac{\lambda}{x_1}$, which implies that $a_0 = c_0 x_1^{\lambda}$ for some $c_0 \in F \setminus \{ 0 \}$. Replacing $g$ with $\frac{g}{c_0}$, which is still a solution of the partial differential equation, we may assume that $c_0 = 1$, so $a_0 = x_1^{\lambda}$.

For $a_1$, we have the equation:
\[(\mu + x_1)a_1 + x_1 a_1' + x_1a_0' = \lambda a_1\]
\noindent which is rewritten as:
\[x_1 a_1' + (x_1 + \mu - \lambda)a_1 + \lambda x_1^{\lambda} = 0 \text{ .}\]

We prove that if this equation has a rational solution, then $\mu$ must be an integer (temporarily replacing $x_1$ with $x$ for less cumbersome notation).

\begin{lemma}\label{lem: LV no x1 Gm}
    If the equation $x a_1' + (x + \mu - \lambda)a_1 + \lambda x^{\lambda} = 0$ has a rational solution, then $\mu \in \mathbb{N}$. 
\end{lemma}

\begin{proof}

    Let $f$ be a rational solution, note that $f \neq 0$ as $\lambda \neq 0$. We see that $f$ has no pole, except maybe at $0$. Therefore its denominator is of the form $x^m$, for some $m \in \mathbb{Z}$, and we write $f = \frac{P}{x^m}$, for some $P \in F[x]$. We have that $P = \sum\limits_{i=0}^n c_i x^i$ with $c_i \in F$ and $c_n \neq 0$ because $f \neq 0$.

    We insert this into the differential equation and obtain:
    \[x\frac{P'x^m - m x^{m-1}P}{x^{2m}} + (x+\mu - \lambda)\frac{P}{x^m} + \lambda x^{\lambda} = 0\]
    \noindent which gives, after some simplifications:
    \[xP'-mP + (x+ \mu - \lambda) P + \lambda x^{\lambda + m} = 0 \text{ .}\]
    This gives us, for all $i \neq \lambda + m$ (with $c_{n+1} = c_{-1} = 0$):
    \[ (i+ \mu -\lambda - m )c_i + c_{i-1} = 0\]
    \noindent which implies that $n = \lambda + m - 1$, as otherwise we would have $c_n = 0$ by picking $i = n+1$.

    The equations for the $c_i$'s imply that for all $i$, if $c_{i-1} = 0$ then $c_i = 0$, as long as $i + \mu - \lambda - m \neq 0$. As $c_n \neq 0$, we must have $\mu = m + \lambda - i$ for some $i \in \{ 0, \cdots , n \}$, which implies $\mu \in \mathbb{N}$.

\end{proof}

To finish the proof, we will need to also consider $g$ as an element of $F(x_0)((x_1))$, meaning $g = \sum\limits_{i = M}^{\infty} b_i x_1^i$, with $b_i \in F(x_0)$ and $b_M \neq 0$. By again inserting this into the differential equation for $g$, we obtain, for all $i$:
\[\mu x_0 b_i' + \left( i(1+x_0) - \lambda \right)b_i + x_0 b_{i-1}'=0\]
\noindent where the derivatives are taken with respect to $x_0$ this time. For $i = M$, this gives us:
\[\frac{b_M'}{b_M} = \frac{\frac{\lambda - M}{\mu}}{x_0} - \frac{M}{\mu}\]
\noindent which implies $M=0$ and $\frac{\lambda}{\mu} \in \mathbb{Z}$. Therefore $b_0 = \alpha x_0^{\frac{\lambda}{\mu}}$ for some $\alpha \in F$.

For $b_1$, we obtain the equation:
\[\mu x_0 b_1 ' + (x_0+1- \lambda)b_1 + \frac{\alpha \lambda}{\mu} x_0^{\frac{\lambda}{\mu}} = 0\]
\noindent and we now prove that this can only have solution if $\frac{1}{\mu} \in \mathbb{N}$ (again temporarily letting $x_0 = x$):

\begin{lemma}\label{lem: LV no x0 Gm}
    If the equation $\mu x b_1 ' + (x+1- \lambda)b_1 + \frac{\alpha \lambda}{\mu} x^{\frac{\lambda}{\mu}} = 0$ has a rational solution, then $\frac{1}{\mu} \in \mathbb{N}$. 
\end{lemma}

\begin{proof}
    Let $f$ be a rational solution. Again $f$ cannot have any pole, expect maybe at $0$. Its denominator must be of the form $x^m$ for some $m \in \mathbb{Z}$, and we can write $f = \frac{Q}{x^m}$, with $Q = \sum\limits_{i = 0}^n d_i x^i$ for some $d_i \in F$ and $d_n \neq 0$.

    We insert this expression into the differential equation and obtain:
    \[\mu x \frac{Q' x^m - m Q x^{m-1}}{x^{2m}} + (x+1 - \lambda) \frac{Q}{x^m} + \frac{\alpha \lambda}{\mu} x^{\frac{\lambda}{\mu}} = 0\]
    \noindent which gives, after simplifications:
    \[\mu x Q' -\mu mQ + (x+1-\lambda) Q + \frac{\alpha \lambda}{\mu} x^{\frac{\lambda}{\mu} + m } = 0 \text{ .}\]
    \noindent We obtain, for all $i \neq \frac{\lambda}{\mu} + m$:
    \[ \left(\mu (i-m) + 1-\lambda \right)d_i + d_{i-1}  = 0\]
    \noindent which implies $n+1 = \frac{\lambda}{\mu} + m$ as otherwise $d_n = 0$ by picking $i = n+1$.

    The equations for the $d_i$ imply that for all $i$, if $d_{i-1} = 0$ then $d_i = 0$, unless $\mu (i-m) + 1-\lambda = 0$. Therefore, for $f$ to be non-zero, we must have $\mu(i-m) + 1 - \lambda =0$ for some $i \in \{ 0, \cdots  ,n\}$, which implies that:
    \begin{align*}
        0 & = i\mu-m\mu + 1 - \lambda \\
        & = i \mu - \mu(n+1) + \lambda + 1 - \lambda \\
        & = -\mu(n+1 - i) + 1
    \end{align*}
    \noindent and thus $\frac{1}{\mu} = n+1 - i \in \mathbb{N}$.
\end{proof}

Combining Lemmas \ref{lem: LV no x1 Gm} and \ref{lem: LV no x0 Gm}, we obtain:

\begin{corollary}\label{cor: LV no Gm}
    There is no solution $g \in F(x_0,x_1)$ to the differential equation $\partials{g}{x_0}\left(\mu x_0 + x_0 x_1\right) +  \partials{g}{x_1} \left( x_1 + x_0 x_1 \right) = \lambda g$ unless $\lambda \in \mathbb{Z}$ and $\mu = 1$. In that case, solutions are a one dimensional $F$-vector space generated by $g = (x_0 - x_1)^{\lambda}$.
\end{corollary}

\begin{proof}
    If there is a solution $g$, then there are solutions to the equations of Lemmas \ref{lem: LV no x1 Gm} and \ref{lem: LV no x0 Gm}, which imply $\lambda \in \mathbb{Z}$ and both $\mu \in \mathbb{N}$ and $\frac{1}{\mu} \in \mathbb{N}$, so $\mu = 1$.

    Conversely, if $\lambda \in \mathbb{Z}$ and $\mu = 1$, it is easy to check that $g = \alpha (x_0 - x_1)^{\lambda}$ is a solution for any $\alpha \in F$. 
    
    Let $h$ be another solution. Recall that this means that $g$ and $h$ induce maps $(u,v) \rightarrow g(u,v)$ and $(u,v) \rightarrow h(u,v)$ from $p$ to the generic type of $z' = \lambda z$. In particular, this implies that $\frac{g(u,v)}{h(u,v)} \in \calC$ for any $(u,v) \models p$. But Lemma \ref{lem: LV weak ortho} tells us that $p$ is weakly $\calC$-orthogonal, so $\frac{g}{f} \in F$.
\end{proof}

We can finally show that the generic type of a Lotka-Volterra system is almost always orthogonal to the constants:

\begin{theorem}\label{theo: LV-is-ortho}
    Let $a,b,c,d \in \mathcal{C} \setminus \{ 0 \}$ and $F = \mathbb{Q}(a,b,c,d)^{\mathrm{alg}}$. The generic type $p$ of the Lotka-Volterra system:
    \[
    \begin{cases}
        x' = ax + b x y \\
        y' = cy + d x y
        \end{cases}
    \]
    is orthogonal to $\calC$ unless $a = c$. In that case, it is 2-analyzable over $\mathcal{C}$, but not almost $\calC$-internal.
\end{theorem}

\begin{proof}
    If $a \neq c$, this is obtained by combining Corollary \ref{cor: ortho crit} with Lemmas \ref{lem: LV weak ortho}, \ref{lem: LV no Ga} and Corollary \ref{cor: LV no Gm}. 

    We now assume that $a = c$. Let $z = \frac{x}{b} - \frac{y}{d}$, it is easy to check that $z' = az$. Thus $(x,y) \rightarrow \frac{x}{b} - \frac{y}{d}$ maps $p$ to solutions of $u' = au$, which is a definable set internal to $\calC$. 

    Given any $(x,y) \models p$ and $z = \frac{x}{b} - \frac{y}{d} $, we can then write:
    \[\begin{cases}
        x = b z + \frac{by}{d}\\
        y' = by^2 + \left( c + db z \right)y
    \end{cases}\]
    \noindent and thus $\tp(x,y/z,F)$ is interdefinable with the generic type of a Bernoulli equation, which is $\calC$-internal. Therefore $p$ is $\calC$-analysable in two steps. To see that $p$ is not $\calC$-internal, note that Lemmas \ref{lem: LV weak ortho} and \ref{lem: LV no Ga} are true with the $a = c$ assumption. Moreover, Corollary \ref{cor: LV no Gm} tells us that if $c\partials{g}{x_0}\left(\mu x_0 + x_0 x_1\right) +  c\partials{g}{x_1} \left( x_1 + x_0 x_1 \right) = \lambda g$ has solutions $g, \lambda$, then $\lambda = k c$ for some $k \in \mathbb{Z}$. In particular, all solutions are $\mathbb{Q}$-linearly dependent. We conclude by Theorem \ref{theo: inter-crit-not-weak-ortho} that $p$ is not almost $\calC$-internal, as there is only one of the two maps needed.
\end{proof}

Some efforts have been made to give exact solutions of Lotka-Volterra systems, see for example \cite{evans1999new}. The fact that the generic solution is 2-analysable in the constants if $a=c$ can be deduced from an article of Varma \cite{varma1977exact}, which in fact shows that it is elementary. More precisely, if $a=c$ Varma proves that the set of solutions is:
\[\begin{cases}
    x = \frac{\alpha e^{at}}{d-be^{\frac{\alpha e^{at}- \beta}{a}}} \\
    y = \frac{\alpha e^{at}}{d e^{-\frac{\alpha e^{at}- \beta}{a}} - b}
\end{cases}\]
\noindent where $\alpha, \beta$ are arbitrary constants.

We were unable to find a proof that solutions cannot be written in terms of elementary functions when $a \neq c$, or a precise statement of what this should mean. We therefore take the opportunity to record that the generic solution, if $a \neq c$, is not classical in the sense of Umemura \cite{umemura1988irreducibility}, meaning, roughly, that it cannot be obtained by solving successive algebraic and logarithmic-differential equations. See the appendix of \cite{nagloo2017algebraic} for more precision. In particular, it is not Liouvillian (see \cite[Definition 4.2]{freitag2023equations}).

\begin{corollary}\label{cor: LV-not-classical}
    If $a,b,c,d \in \calC \setminus \{ 0 \}$ then the generic solution of the Lotka-Volterra system:
    \[
    \begin{cases}
        x' = ax + b x y \\
        y' = cy + d x y
        \end{cases}
    \]
    \noindent is classical if and only if $a = c$, in which case it is elementary.
\end{corollary}

\begin{proof}
    If $a \neq c$, then by Theorem \ref{theo: LV-is-ortho}, the generic type of the Lotka-Volterra system is orthogonal to the constants, and therefore not $\calC$-analysable. It follows by \cite[Proposition A.3]{nagloo2017algebraic} that the generic solution is not classical.

    Now suppose that $a = c$, then by \cite{varma1977exact}, the solutions are elementary.

\end{proof}

Since this article was first written, Duan and Nagloo \cite{duan2025algebraic} obtained a much stronger result if $\frac{a}{c} \not\in \mathbb{Q}$: the type $p$ is then \emph{strongly minimal}. This was then improved to $a \neq c$ in \cite{duan2025algebraicm} by Duan and the authors. However, the computations and results of this subsection are still necessary to obtain the complete classification of invariant algebraic curves of \cite[Theorem 3.13]{duan2025algebraicm}, in the $a = c$ case.

\subsection{Pullbacks}

It could be expected that the work of this article subsumes some of the results of \cite{jin2020internality} and \cite{eagles2024splitting} on internality of logarithmic-differential pullbacks. In this subsection, we briefly explain why this is not the case. 

Both articles consider the generic type $p$ of some differential equation: 
\[y^{(n)} = f(y, y' , \cdots , y^{(n-1)})\]
\noindent where $f \in F(x)$ and $F$ is some algebraically closed differential field. One can consider the generic type $q : = \logd^{-1}(p)$ of the system of equations:
\[
\begin{cases}
    y^{(n)} = f(y, y' , \cdots , y^{(n-1)}) \\
    z' = yz
\end{cases}
\]
\noindent which is the generic type of the pullback of the first equation under the logarithmic derivative. The question these articles ask is: assuming that $p$ is almost $\calC$-internal, when is $q$ also almost $\calC$-internal?

Over constant parameters, the answer given is that $q$ is almost $\calC$-internal if and only if there are a non-zero $h \in F(x_0, \cdots , x_{n-1})$, some $e \in F$ and some integer $k \neq 0$ such that:

\begin{align*} (kx_0-e)h = \sum\limits_{i=0}^{n-2} \partials{h}{x_i} x_{i+1} + \partials{h}{x_{n-1}}f \tag{$\ast$} \label{ast} \text{ .}\end{align*}

Using Theorem \ref{theo: intern-crit-weak-ortho} directly can only give the existence of some non-zero $g \in F(x_0, \cdots , x_{n-1},x_n)$ and $\lambda \in F$ such that:

\[\lambda g = \sum\limits_{i=0}^{n-2} \partials{g}{x_i} x_{i+1} + \partials{g}{x_{n-1}}f + \partials{g}{x_n}x_0x_n \]

\noindent which has an extra variable. The only extra information we could use to eliminate $x_n$ is that the map given by $g$ does not factor through $p$, and thus $\partials{g}{x_n} \neq 0$. This does not seem to be enough to recover (\ref{ast}). What is missing is the following: in \cite{eagles2024splitting} it is proven (under various extra assumptions, one of them being that the equation is autonomous) that almost $\calC$-internality of $q$ is equivalent to the existence of an integer $k \neq 0$ such that for some $a \models p$, there are $w_1,w_2$ with:
\begin{itemize}
    \item $a^k = w_1 w_2$
    \item $w_1 \in \dcl(F, \logd a)$
    \item $\logd (w_2) \in F$.
\end{itemize}

\noindent This is key to obtaining (\ref{ast}), and does not seem to follow easily from Theorem \ref{theo: intern-crit-weak-ortho}.

\section{Remarks on generalizations}\label{sec: generalizations}

There are two natural potential generalizations of our work:
\begin{enumerate}[(a)]
    \item instead of considering rational vector fields on affine space, work with rational vector fields on an arbitrary algebraic variety,
    \item consider non-autonomous systems.
\end{enumerate}

The main obstacle to considering differential systems on arbitrary varieties is that we cannot guarantee that the binding group is linear. Looking at the proof of Lemma \ref{lem: lin-bind-group-im}, it should be possible to obtain similar results for equations on unirational algebraic varieties. 

Beyond that, the binding group may not be linear, but should always have a commutative linear part. It is possible that our methods would then generalize.

For non-autonomous systems, the main obstacle is also Galois-theoretic. Indeed, we completely loose control over the binding group in that case: it is a theorem of Kolchin (see \cite[Theorem 2]{kolchin1955galois}) that any algebraic group can occur as a binding group. 

For a single order one equation, Jaoui and Moosa proved (see \cite[Proposition 6.5]{jaoui2022abelian}): 

\begin{theorem}[Jaoui-Moosa]
    Let $F$ be an algebraically closed differential field. Let $p \in S(F)$ be one-dimensional, $\calC$-internal, weakly $\calC$-orthogonal and with the binding group $\bg{F}{p}$ equal to its linear part. Then $p$ is interdefinable with the generic type of a Ricatti equation $y' = ay^2 + by + c$.
\end{theorem}

At the heart of their proof is the fact that the type $p$ must be strongly minimal, and thus by a result of Hrushovski \cite{hrushovski1989almost}, the action of $\bg{F}{p}$ on $p$ must be definably isomorphic to either the action of $G_a(\calC), G_m(\calC)$ or the constant points of an elliptic curve on themselves, the natural action of $\mathrm{Aff}_2(\calC)$ on $\calC$, or the natural action of $\mathrm{PSL}_2(\calC)$ on $\mathbb{P}^1(\calC)$. If the binding group is linear, it rules out an elliptic curve, and they reduce all other cases to a Ricatti equation. 

We record the following corollary, which was not observed in \cite{jaoui2022abelian}. It gives a Rosenlicht style criteria for almost $\calC$-internality, and indeed recovers Rosenlicht's theorem \cite{rosenlicht1974nonminimality} as special cases:

\begin{corollary}\label{cor: Jaoui-Moosa-cor}
    Let $F$ be an algebraically closed differential field, some $f \in F(x)$, and $p$ the generic type of the equation $y ' = f(y)$. Then $p$ is almost $\calC$-internal if and only if there is a non-constant $g \in F(x)$ and some $a,b,c \in F$ such that:
    \[f \partials{g}{x} + g^{\delta}= a g^2 + bg + c \text{ .}\]
\end{corollary}

\begin{proof}
    First note that we may assume that $p$ is weakly $\calC$-orthogonal as $p$ being not weakly $\calC$-orthogonal corresponds to the case $a=b=c = 0$.
    The proof is then entirely similar to that of Theorem \ref{theo: intern-crit-weak-ortho}, but using the map given by Jaoui and Moosa's theorem.
    
    For the converse, note that $g$ induces an $F$-definable map from $p$ to the generic type of the Ricatti equation $z' = az^2 + bz + c$, which is $\calC$-internal and of $U$-rank one. Thus the induced map is finite-to-one, and $p$ is almost $\calC$-internal.
\end{proof}

Another way to interpret Jaoui and Moosa's proof is as follows: the action of $\bg{F}{p}$ is definably isomorphic to a group of birational transformations of $\mathbb{P}^1$. But the group of birational transformations of $\mathbb{P}^1$ is equal to $\mathrm{PGL}_2$, and their proof examines every possible connected algebraic subgroup of $\mathrm{PGL}_2$ (we know the binding group is connected because the parameters are algebraically closed). In higher dimension, one would need to look at the non-algebraic Cremona groups and their possible connected algebraic subgroups. A classification of the maximal connected algebraic subgroups has been obtained in dimension two by Enriques \cite{enriques1893sui} (see also \cite{umemura1982maximal}), and it is simple enough that a complete characterization of internality for non-autonomous systems of differential equations on the plane may be obtained. We reserve this for future work.

\medskip

\textbf{Acknowledgements.} The authors are thankful to Rahim Moosa for multiple conversations on the subject of this article. The second author is also grateful to R\'{e}mi Jaoui for a related discussion. Some of this work was done while the second author was visiting the University of Waterloo: he is thankful to Nicolas Chavarria Gomez for hosting him in his office. The authors are also grateful to Gleb Pogudin for suggesting an improvement to the main theorem by using Lie derivatives. Finally, we thank the anonymous referees for suggesting numerous improvements.

\bibliography{biblio}
\bibliographystyle{plain}

\end{document}